\documentclass[pdflatex,sn-mathphys-num]{sn-jnl}

\usepackage{graphicx}%
\usepackage{multirow}%
\usepackage{amsmath,amssymb,amsfonts}%
\usepackage{amsthm}%
\usepackage{mathrsfs}%
\usepackage[title]{appendix}%
\usepackage{xcolor}%
\usepackage{textcomp}%
\usepackage{manyfoot}%
\usepackage{booktabs}%
\usepackage{algorithm}%
\usepackage{algorithmicx}%
\usepackage{algpseudocode}%
\usepackage{listings}%

\usepackage[caption=false]{subfig}

\theoremstyle{thmstyleone}%
\newtheorem{theorem}{Theorem}
\newtheorem{lemma}[theorem]{Lemma}%
\newtheorem{corollary}[theorem]{Corollary}%

\theoremstyle{thmstyletwo}%
\newtheorem{example}{Example}%
\newtheorem{remark}{Remark}%

\theoremstyle{thmstylethree}%

\newcommand{\bomega}{\boldsymbol \omega}

\begin{document}

\title[Random Batch Method for Networked 1-D Wave Equations]{A Random Batch Method for the Efficient Simulation and Optimal Control of Networked 1-D Wave Equations}


\author*[1]{\fnm{Dani\"el W.~M.} \sur{Veldman}}\email{d.w.m.veldman@gmail.com}

\author*[2]{\fnm{Yue} \sur{Wang}}\email{yuewang@fudan.edu.cn}
\equalcont{These authors contributed equally to this work.}

\affil[1]{\orgdiv{Department of Mathematics}, \orgname{Friedrich Alexander Universit\"at Erlangen-Neuremberg}, \orgaddress{\street{Cauerstrasse 11}, \city{Erlangen}, \postcode{91052}, \state{Bavaria}, \country{Germany}}}

\affil[2]{\orgdiv{Center for Applied Mathematics \& Department of Mathematics}, \orgname{Fudan University}, \orgaddress{\street{Handan Road 220}, \city{Shanghai}, \postcode{200433}, \state{Shanghai}, \country{China}}}


\abstract{In this paper, a stochastic algorithm for the efficient simulation and optimal control of networked wave equations based on the random batch method is proposed and analyzed. The random approximation is constructed by dividing the time interval into subintervals and restricting the dynamics to a randomly chosen subnetwork during each of these subintervals. It is proven that the solution for this randomized system converges in expectation to the solution on the original network when the length of the subintervals approaches zero. Furthermore, the optimal controls for the randomized system converge (in $H^2$ and in expectation) to the optimal controls for the original system.  
The computational advantage of the proposed method is demonstrated in two numerical examples. }

\keywords{Random Batch Method, Stochastic Algorithms, Wave Equation, Efficient Simulation, Optimal Control, Domain Decomposition}



\maketitle

\section{Introduction}
Networks of hyperbolic systems can describe a wide range of physical phenomena, e.g., Euler equations of gas dynamics, shallow water equations, Maxwell equations, magnetohydrodynamics equations, and classical and elastic wave equations. In many applications, models consisting of networked one-dimensional hyperbolic systems are appropriate and these models often also need to be controlled.  Some recent examples include the network of vibrating strings \cite{YueLi2021nodal, Leugering2019LW}, the network of gas pipe lines in Germany \cite{herty2007}, networks of geometrically exact beams that model a flexible aircraft \cite{leugering_beams2021}, and the control of the water level in a canal with an immersed obstacle \cite{Vergara2021}. 
Especially when the considered networks are large, the computational cost for simulation or optimal control can be large, see e.g. \cite{hante2017}. Developing numerically efficient algorithms for the simulation and control of such networks is therefore an important task. 

Inspired by the successes of stochastic algorithms such as Stochastic Gradient Descent (SGD) in supervised learning, Random Batch Methods (RBMs) have been introduced recently in \cite{jin2020} for the efficient simulation of interacting particle systems. Because there are $N(N-1)/2$ interactions between $N$ particles, the computational cost for the simulation of an $N$-particle system is of $O(N^2)$ and thus grows rapidly with $N$. By considering not all $N-1$ interactions a particle has with all the other particles, but only $p \ll N$ randomly chosen interactions, the RBM reduces the computational cost from $O(N^2)$ to $O(pN)$. When the interactions are chosen independently in each time interval, it can be shown that the RBM-solution converges to the solution of the original problem (that considers all interactions between particles) when the length of the time intervals in which the randomization is refreshed approaches zero. 

As demonstrated in \cite{ko2021model}, the same ideas can also be used speed up the solution of optimal control problems governed by interacting particles systems. Because the number of interactions considered in the adjoint system also reduces from $O(N^2)$ to $O(pN)$, the RBM also significantly the computational cost of each iteration in a gradient-based algorithm. The results in \cite{ko2021model} demonstrate that a combination of Model Predictive Control (MPC) with the RBM can lead to the accurate and fast approximation of the optimal control for large-scale interacting particles systems. However, the convergence of the RBM-optimal control to the optimal control for the original system was not proven in \cite{ko2021model}. 

A first rigorous analysis of RBMs in linear-quadratic optimal control was provided in \cite{veldman2022}. Here it was proven that the RBM-optimal control converges to the control of the original system when the length of the time intervals approaches zero. The same paper also shows that original idea of the RBM for interacting particle systems can be extended to a more general operator-splitting setting. Although numerical simulations in \cite{veldman2022} suggested that the RBM is also applicable to infinite dimensional systems, the analysis in \cite{veldman2022} was limited to finite-dimensional systems. Recently, it has been shown in \cite{eisenmann2022} that the RBM is also applicable to a backward Euler discretization of a certain class of nonlinear parabolic Partial Differential Equations (PDEs) that for example includes the $p$-Laplacian. The results in \cite{eisenmann2022} only consider the forward dynamics and do not consider optimal control. In contrast, the RBM for parabolic PDEs on graphs was studied in \cite{hernandez2025} by considering a spatial discretization but keeping time continuous. Convergence of the RBM dynamics and the corresponding optimal control problem can then be proved by invoking classical convergence results for the spatial discretization and the convergence results for the RBM in finite dimensions from \cite{veldman2022}. 

This paper is the first application of the RBM to hyperbolic PDEs. As a prototypical example a network of (classical) wave equations is considered. The RBM-approximation is constructed as follows. First, the considered time interval is divided into small sub-intervals. In each subinterval, a subset of edges of the network is randomly chosen. The velocity of propagation on the other edges is set to zero, whereas the velocity of propagation on the chosen edges is set such that expected value of the velocity of propagation is equal to the velocity of propagation in the original network. We prove that the thus constructed RBM-approximation converges to the solution on the original network when the length of the sub-intervals $h$ approaches zero. Furthermore, it is proven that the optimal controls for the ranomized network converge to the optimal controls for the original system (in expectation and in $H^2$) for $h \rightarrow 0$. The effectiveness of the proposed method is demonstrated in two numerical examples. 

In contrast to previously developed methods for infinite-dimensional systems \cite{eisenmann2022, hernandez2025}, our analysis does not rely on the discretization of space or time.

The proposed algorithm can also be viewed as a randomized version of nonoverlapping domain decomposition for hyperbolic systems, see, for example, \cite{nataf2015}. Randomized domain decomposition for elliptic problems has been studied explicitly in \cite{griebel2012} and for problems of parabolic type in \cite{eisenmann2022}. The application of domain decomposition to optimal control problems is now well-developed and can be traced back to Glowinski and J.L. Lions, see, for example, \cite{Bensoussan1973}. Later, new techniques for space-time decomposition of the optimality system have been proposed by Lagnese and Leugering \cite{lagnese2003, lagnese2012}. These techniques have also been applied to wave equations, see, for example, \cite{krug2021}, and networks of 1D elements, see, for example, \cite{leugering2023}. Recently, domain decomposition for optimal control problems also received attention from other authors, see, for example, \cite{gander2024, cocquet2025}. We are not aware of any existing publications on randomized domain decomposition for hyperbolic PDEs and the assiciated optimal control problems. 

The remainder of this paper is structured as follows. In Section \ref{sec:problem_results}, the RBM-approximation is introduced in more detail and the convergence results are presented. The convergence for the RBM-dynamics is proven in Section \ref{sec:dynamics} and the convergence in optimal controls is presented in Section \ref{sec:control}. Section \ref{sec:examples} contains two numerical examples that demonstrate the effectiveness of the proposed method and Section \ref{sec:conclusions} contains conclusions and an outlook to future work. 

\section{Problem Setting and Main Results \label{sec:problem_results}}

\subsection{Modeling of Networked Wave Equations and Optimal Control}

In this section, the considered model for networked linear wave equations on directed metric graphs is described. The graph $G= (V, E, L)$ consists of a set of vertices $V$ which are enumerated as $v_j$ ($j \in \{1,2,\ldots,|V| \}$), a set of edges $E$ enumerated as $e_i$ ($i \in \{1,2, \ldots,|E| \}$), and a tuple $L = \mathbb R ^ {|E|}$ of lengths $\ell_{e_i}$ for each edge $e_i$. The set of edges connected to the vertex $v_j \in V$ is denoted by $E(v_j)$ and the degree of $v_j$ is $|E(v_j)|$.  
The set of edges is then also described by the incidence matrix $D \in \mathbb R^{|V| \times |E|}$ with the element defined as
\begin{equation} \label{eq.index_matrix}
    D_{ji}  = \left \{ 
    \begin{aligned} 
     -1, \quad & \text{if node $v_j$ is the start point of edge $e_i$;} \\
     1, \quad & \text{if node $v_j$ is the end point of edge $e_i$;} \\
     0, \quad &\text{otherwise},
    \end{aligned}
    \right.
\end{equation}

\begin{example} \label{example:diamond_graph}(Diamond directed Graph) Consider the diamond network in Figure \ref{fig:diamond-graph} with vertices $V = \{v_1, v_2, ..., v_6\}$ and edges $E = \{e_1,...,e_7\}$.  
With the directions for the edges as in Figure \ref{fig:diamond-graph}, the incidence matrix becomes
\begin{equation}
    D =  \begin{bmatrix}
            -1 & 0 & 0 & 0 & 0 & 0 & 0\\
            1 & -1 & -1 & 0 & 0 & 0 & 0\\
            0 & 1 & 0 & -1 & -1 & 0 & 0\\
            0 & 0 & 1 & 1 & 0 & -1 & 0\\
            0 & 0 & 0 & 0 & 1 & 1 & -1\\
            0 & 0 & 0 & 0 & 0 & 0 & 1\\
    \end{bmatrix}.
\end{equation}

\begin{figure}[h]
    \centering
    \includegraphics[width=0.7
    \textwidth]{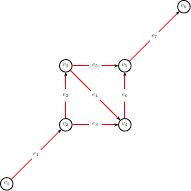}
    \caption{The diamond graph considered in Example \ref{example:diamond_graph}}
    \label{fig:diamond-graph}
\end{figure}
\end{example}

On each edge $e_i$, a coordinate $x \in [0, \ell_{e_i}]$ is introduced for which $x = 0$ corresponds to the starting point and $x = \ell_{e_i}$ corresponds to the end point. The transversal displacement of the string on the edge $e_i \in E$ is denoted $y^{e_i} : [0,\infty) \times [0, \ell_{e_i}] \mapsto \mathbb R$. The solution $\mathbf{y} = (y^{e_i}(t,x))_{e_i \in E}$ satisfies the following equations
\begin{equation}\label{eq.system}
\left\{\begin{aligned}
    &y^{e_i}_{tt} (t,x) - c_{e_i}^2 y^{e_i}_{xx} (t,x)=0, \quad & e_i\in E, \\
    &\sum_{e_i\in E(v_j)} c_{e_i}^2 D_{ji} y^{e_i}_x (t, v_j) =  \bar{u}^{v_j}(t),  & v_j \in V, \\
    & y^{e_i} (t, v_j) = y^{e_k} (t, v_j),\quad & \forall e_i,e_k \in E(v_j), v_j \in V,\\
    & y^{e_i} (0,x) = y_0^{e_i}(x), \quad y_t^{e_i} (0,x) = y^{e_i}_1 (x), \quad &  e_i\in E,  
\end{aligned} 
\right.
\end{equation}
where $c_{e_i} > 0$ is the velocity of propagation on edge $e_i$, 
\begin{equation}
    \bar{u}^{v_j}(t) = \begin{cases}
        - u^{v_j} (t), \quad  & v_j \in V_C, \\
        0,  \quad & v_j \in V\backslash V_C,
    \end{cases}
\end{equation}
is the external Neumann (force) control $\mathbf{u} = (u^{v_j}(t))_{v_j \in V_C}$ acting only on vertices $v_j \in V_C \subseteq V$, and $\mathbf{y}_0 = (y_0^{e_i}(x))_{e_i \in E}$ and $\mathbf{y}_1 = (y_1^{e_i}(x))_{e_i \in E}$ are the initial conditions. The conditions in the second line are of Kirchoff-type and the conditions in the third line assure continuity at the nodes (geometric multiple node condition, see, e.g.,  \cite{lagnese2012modeling}). It is required that $y_{0,x}^{e_i}$ and $y_{1}^{e_i}$ are uniformly Lipschitz continuous on each edge $e_i \in E$ and satisfy standard compatibility conditions at the boundaries. 

Let $T>0$, $\Omega = \Pi_{e_i\in E} (0, l_{e_i})$, and $Q = (0,T) \times \Omega$. For given functions $\mathbf y_d \in L^2 (Q)$ and a weight $\alpha > 0$, we also consider the following optimal control problem 
\begin{equation}
    \min_{\mathbf u \in U_{ad}} J (\mathbf u) = \frac{1}{2} \| \mathbf y - \mathbf y_d \|^2_{L^2(Q)} + \frac{\alpha}{2} |\mathbf{u}|_{H^2(0,T)}^2 \label{eq:OCP} \tag{OCP}
\end{equation}
where $U_{ad}$ is a convex closed subset of $H^2(0,T;\mathbb{R}^{|V_C|})$,
$\mathbf y$ is the solution of \eqref{eq.system} 
resulting from the initial conditions $\mathbf y_0$ and $\mathbf y_1$ and the control $\mathbf u$, which can also be obtained by solving \eqref{eq:D} and using \eqref{eq:yintw}. 
The minimizer of \eqref{eq:OCP} is denoted by $\mathbf u^*$.

\subsection{Riemann invariants}

In order to prepare for the introduction of the randomized dynamics in the following subsection, we transform the second order dynamics in \eqref{eq.system} into a first-order system in terms of the Riemann invariants
\begin{equation}\label{eq.Riemann}
  \mathbf{w}^{e_i}(t,x) = \begin{pmatrix}
   w^{e_i}_-(t,x) \\ w^{e_i}_+(t,x)   
  \end{pmatrix} = \begin{pmatrix}
    y^{e_i}_t(t,x) + c_{e_i} y^{e_i}_x(t,x)  \\
    y^{e_i}_t(t,x) - c_{e_i} y^{e_i}_x(t,x)
\end{pmatrix}, 
\end{equation}
which satisfy
\begin{equation} \left\{
\begin{aligned}
& w_{-,t}^{e_i}(t,x) - c_{e_i} w_{-,x}^{e_i}(t,x) = 0, \\
& w_{+,t}^{e_i}(t,x) + c_{e_i} w_{+,x}^{e_i}(t,x) = 0.
\end{aligned}\right.
\end{equation}


The boundary conditions can be rewritten by the same transformation. One easily verifies that the force condition at each vertex $v_j\in V$ becomes
\begin{equation}\label{eq.BC_flux}
    \sum_{e_i\in E(v_j)} c_{e_i} D_{ji} \frac{w^{e_i}_-(t, v_j) - w^{e_i}_+(t, v_j)}{2} = \bar{u}^{v_j}(t)
\end{equation}
and the continuity constraint at the junction becomes
\begin{equation}\label{eq.BC_contin}
    w^{e_i}_-(t, v_j) + w^{e_i}_+(t, v_j) = w^{e_k}_-(t, v_j) + w^{e_k}_+(t, v_j), \quad  e_i,e_k\in E(v_j), i\neq k.
\end{equation}

Furthermore, for each $e_i\in E(v_j)$, the value of Riemann invariant leaving the edge $e_i$ at the vertex $v_j$ is 
\begin{subequations} \label{eq:riemann_inout}
    \begin{align}
        &w_{\mathrm{out}}^{e_i}(t, v_j) = y^{e_i}_t (t, v_j) + D_{ji} c_{e_i} y^{e_i}_x (t, v_j), 
        \end{align}
and the value of the Riemann invariant entering the edge $e_i$ at node $v_j$ is
\begin{align}
&w_{\mathrm{in}}^{e_i}(t, v_j) = y^{e_i}_t (t, v_j) - D_{ji} c_{e_i} y^{e_i}_x (t, v_j).
    \end{align}
\end{subequations} Hence, combining the boundary/transmission conditions \eqref{eq.BC_flux} and \eqref{eq.BC_contin}, each `entering-flow' $w_{\mathrm{in}}^{e_i}(t,v_j)
$ can be expressed in terms of all `leaving-flows' $w_{\mathrm{out}}^{e_k}(t,v_j)$ at the node $v_j$ ($e_k \in E(v_j)$). To be specific, \eqref{eq.BC_flux} and \eqref{eq.BC_contin} can be rewritten as (see Appendix \ref{app:bcs})
\begin{equation}
w^{e_i}_{\mathrm{in}}(t,v_j) = - w^{e_i}_{\mathrm{out}}(t,v_j) + \frac{2}{c^{v_j}_{\mathrm{tot}}}\left( \sum_{e_k \in E(v_j)}  c_{e_k}w_{\mathrm{out},j}^{e_k}(t,v_j) - \bar{u}^{v_j}(t) \right),
\end{equation}
where
\begin{equation} \label{eq:cvjtot}
c^{v_j}_{\mathrm{tot}} =  \sum_{e_k \in E(v_j)} c_{e_k}.
\end{equation}

Thus, the system \eqref{eq.system} can be rewritten as
\begin{equation} \left\{
\begin{aligned}
& w_{-,t}^{e_i}(t,x) - c_{e_i} w_{-,x}^{e_i}(t,x) = 0, \\
& w_{+,t}^{e_i}(t,x) + c_{e_i} w_{+,x}^{e_i}(t,x) = 0, \\
    & w^{e_i}_{\mathrm{in}}(t,v_j) = - w^{e_i}_{\mathrm{out}}(t,v_j) + \frac{2}{c^{v_j}_{\mathrm{tot}}}\left( \sum_{e_k \in E(v_j)} c_{e_k} w_{\mathrm{out},j}^{e_k}(t,v_j) -\bar{u}^{v_j}(t)\right), \\
    & w^{e_i}_-(0,x) =  y^{e_i}_1 (x) + c_{e_i} y_{0,x}^{e_i}(x), \\ 
    & w^{e_i}_+(0,x) =  y^{e_i}_1 (x) -c_{e_i} y_{0,x}^{e_i}(x)
    \end{aligned}\right. 
    \tag{D} \label{eq:D}
\end{equation}
where the boundary conditions hold for all $e_i \in E(v_j)$. 
Note that the solution $\mathbf{y}$ of \eqref{eq.system} can be obtained by observing that $2y_t^{e_i}(t,x) = w^{e_i}_+(t,x) + w^{e_i}_-(t,x)$, so that
\begin{equation}
    y^{e_i}(t,x) = y_0^{e_i}(x) + \int_0^t \frac{ w^{e_i}_+(s,x) + w^{e_i}_-(s,x)}{2} \ \mathrm{d}s. \label{eq:yintw}
\end{equation}

\subsection{Randomized Dynamics and Optimal Control}
    In order to introduce the randomized approximation, start by enumerating the subsets of $E$ as $E_1, E_2, \ldots , E_{2^{|E|}}$. To each of the subsets $E_\omega$ ($\omega \in \{1,2,\ldots, 2^{|E|} \}$), a probability $p_\omega \geq 0$ is assigned with which this subset will be selected. These probablilties $p_\omega \geq 0$ should satisfy
    \begin{equation}
     \sum_{\omega=1}^{2^M} p_\omega = 1.  \label{eq:sump}
    \end{equation}
    Next, the time interval $[0,T]$ is divided into $K$ subintervals $(t_{k-1}, t_k]$ (with $k \in \{1,2,\ldots, K \}$). The length of these subintervals is bounded by
    \begin{equation}
        h = \max_{k \in \{1, 2, \ldots, K \}} t_k - t_{k-1}. 
    \end{equation}
    For each time interval $(t_{k-1}, t_k]$, an index $\omega_k \in \{1,2,\ldots,2^{|E|} \}$ is selected according to the probabilities $p_\omega$. All indices are collected in a vector
    \begin{equation}
        \bomega = (\omega_1, \omega_2, \ldots, \omega_K) \in \{1,2,\ldots,2^{|E|} \}^K. 
    \end{equation}
    Write
\begin{equation}
        \chi_{e_i}(\omega) =  \begin{cases}  
        1, & e_i \in E_{\omega}, \\
        0, & e_i \notin E_{\omega}.
        \end{cases}, \qquad \qquad \pi_{e_i} := \mathbb{E}[ \chi_{e_i}] = \sum_{\omega \in \{ \omega \mid e_i \in E_{\omega} \}} p_\omega. \label{eq:def_pi}
    \end{equation}
    Note that $\pi_{e_i} \in [0,1]$ represents the probability that an edge $e_i$ is an element of the selected subset. For the construction of the randomized dynamics, it is required that the probabilities $p_\omega$ are chosen such that $\pi_{e_i} > 0$ for all $e_i \in E$. 
    
    Note that each subset $E_\omega$ gives rise to a metric subgraph $G_\omega = (V,E_\omega,L_\omega)$, where $L_\omega = \{\ell_{e_i} \mid e_i \in E_\omega \}$ of the original metric graph $G = (V, E, L)$. The random approximation is constructed by restricting the evolution of the Riemann variables $w_+^{e_i}(t,x)$ and $w^{e_i}_-(t,x)$ to the randomly chosen subgraph $(V,E_{\omega_k},L_{\omega_k})$ in the time interval $(t_{k-1}, t_k]$. When these subgraphs have a much simpler structure than the original graph, this can lead to a significant reduction in computational cost. When $h$ is sufficiently small, i.e. when the switching is sufficiently fast, this will lead to a good approximation of the original dynamics on the metric graph $G = (V, E, L)$.  

    In order to make this idea more precise, introduce
    \begin{equation}
        c_{h,e_i}(\bomega,t) := \frac{c_{e_i}}{\pi_{e_i}} \chi_{e_i}(\omega_k), \qquad \qquad t \in (t_{k-1}, t_k]. 
    \end{equation}
It will be shown in Section \ref{sec:preliminaries} that the expected values of $c_{h,e_i}(\bomega,t)$ is $c_{e_i}$ for all $t \in (0,T]$. 

The randomized approximation for the dynamics in Riemann coordinates in \eqref{eq:D} is then defined for $t \in (t_{k-1}, t_k]$ as
\begin{equation} \left\{
\begin{aligned}
& w_{h-,t}^{e_i}(\bomega,t,x) - c_{h,e_i}(\bomega,t) w_{h-,x}^{e_i}(\bomega,t,x) = 0, \\
& w_{h+,t}^{e_i}(\bomega,t,x) + c_{h,e_i}(\bomega,t) w_{h+,x}^{e_i}(\bomega,t,x) = 0, \\
    & w^{e_i}_{h,\mathrm{in}}(\bomega,t,v_j) = - w^{e_i}_{h,\mathrm{out}}(\bomega,t,v_j) + \frac{2}{c^{v_j}_{\mathrm{tot}}} \left(\sum_{e_k \in E(v_j)} c_{e_k} w^{e_k}_{h,\mathrm{out},j}(\bomega,t,v_j) -  \bar{u}^{v_j}(t) \right), \\
    & w_{h-}^{e_i}(\bomega,0,x) = y^{e_i}_1 (x) + c_{e_i} y_{0,x}^{e_i}(x), \\
    & w_{h+}^{e_i}(\bomega,0,x) = y^{e_i}_1 (x) - c_{e_i} y_{0,x}^{e_i}(x),
    \end{aligned}\right. 
    \tag{RD} \label{eq:RD}
\end{equation}
Note that the outflows $w^{e_k}_{h,\mathrm{out},j}(\bomega,t,v_j)$ of edges $e_k$ for which $c_{h,e_k}(\bomega,t) = 0$ still appear in the boundary conditions and provide the required coupling between the edges/subdomains. 
Just as the solutions $w_{+}^{e_i}(t,x)$ and $w_{-}^{e_i}(t,x)$ determine the solution $\mathbf y$ through \eqref{eq:yintw}, the solutions $w_{h+}^{e_i}(\bomega,t,x)$ and $w_{h-}^{e_i}(\bomega,t,x)$ define the randomized solution $\mathbf{y}_h^{e_i}(\bomega) = (y_h^{e_i}(\bomega,t,x))_{e_i \in E}$ as
\begin{equation}
    y_h^{e_i}(\bomega,t,x) = y_0^{e_i}(x) + \int_0^t \frac{w^{e_i}_{h+}(\bomega,s,x) + w^{e_i}_{h-}(\bomega,s,x)}{2} \ \mathrm{d}s. \label{eq:yhintwh}
\end{equation}

The randomized optimal control problem associated to these dynamics now reads
\begin{equation}
    \min_{\mathbf u \in U_{ad}} J_h(\bomega, \mathbf u) = \frac{1}{2} \| \mathbf y_h(\bomega) - \mathbf y_d \|^2_{L^2(Q)} +\frac{\alpha}{2} |\mathbf{u}|_{H^2(0,T)}^2, \label{eq:OCPh} \tag{ROCP}
\end{equation}
where $U_{ad}$ is as in \eqref{eq:OCP}, 
$\mathbf y_h(\bomega)$ is the solution computed from \eqref{eq:RD} and \eqref{eq:yhintwh}  
resulting from the initial conditions $\mathbf y_0$ and $\mathbf y_1$ and the control $\mathbf u$. 
The minimizer of \eqref{eq:OCPh} is denoted by $\mathbf u^*_h(\bomega,t)$.
Problem \eqref{eq:OCPh} can be much easier to solve than \eqref{eq:OCP} when the randomized forward dynamics \eqref{eq:RD} is much easier to compute than \eqref{eq:D}. 

This subsection is concluded by two examples that illustrate the character of the constructed randomized solutions and the construction of a randomized approximation on more complex network structures. 

\begin{example}[] \label{example:jump}
The simplest situation in which our method can be applied is a graph consisting of three nodes linearly interconnected by two edges. Figure \ref{fig:Jump} now shows the evolution of the randomized solution starting from an harmonic initial condition for the case where the edges $e_1$ and $e_2$ are activated in an alternating manner. From the figure it is clear that although the initial condition is continuous and the the solution to the original dynamics \eqref{eq:D} will therefore also be continuous, the randomized approximation \eqref{eq:RD} introduces discontinuities due to the jump in the velocity at the connecting node. Note that for this very simple network structure without loops, we do not expect to gain a computational advantage of the randomized approximation. However, for more complex network structures as in the following example, the randomized approximation can offer a significant reduction in computational cost and the introduction of discontinuities can be a price worth paying. This will be demonstrated further by the two numerical examples in Section \ref{sec:examples}.

\begin{figure}
\centering

\subfloat[$t = t_0 = 0$ \label{fig:Jumpa}]{
\includegraphics[width=0.45\textwidth]{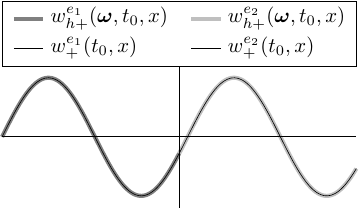}}
~
\subfloat[$t = t_1$ \label{fig:Jumpb}]{
\includegraphics[width=0.45\textwidth]{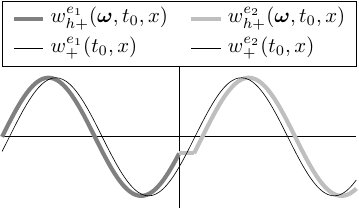}}

\subfloat[$t = t_2$ \label{fig:Jumpc}]{
\includegraphics[width=0.45\textwidth]{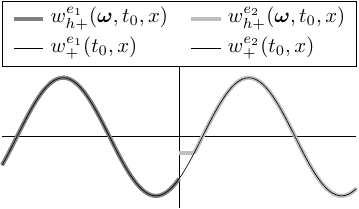}}
~
\subfloat[$t = t_3$ \label{fig:Jumpd}]{
\includegraphics[width=0.45\textwidth]{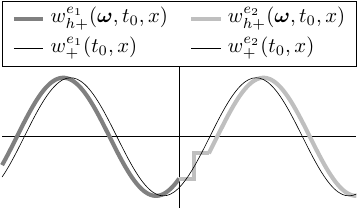}}

\subfloat[$t = t_4$ \label{fig:Jumpe}]{
\includegraphics[width=0.45\textwidth]{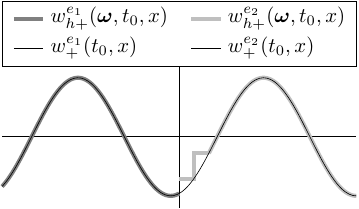}}
~
\subfloat[$t = t_5$ \label{fig:Jumpf}]{
\includegraphics[width=0.45\textwidth]{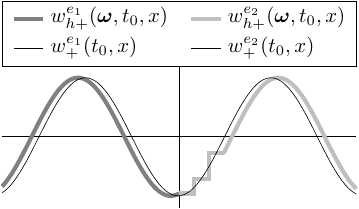}}

\caption{An example of how the activation and disabling of edges typically will introduce discontinuities in the Riemann invariants. Consider the evolution of $w_{h+}^{e_i}(\bomega,t,x)$ on two edges $e_1$ and $e_2$. The initial condition at $t = 0$ is continuous. During $[t_0, t_1]$, edge $e_1$ is disabled and $e_2$ is active and the continuity of $w_{h+}^{e_i}(\bomega,t,x)$ is preserved. During $[t_1, t_2]$, $e_1$ is active and $e_2$ is disabled and a discontinuity at the node connecting the two edges appears. During $[t_2, t_3]$, edge $e_1$ is disabled and $e_2$ is active again and the discontinuity is propagated into the edge $e_2$. During $[t_3, t_5]$ the same process is repeated. The exact solution $w^{e_i}(t,x)$ is also shown for reference. Note that the alternating activation pattern is only chosen as an illustration and that the activation pattern in the proposed strategy is chosen randomly.   }
\label{fig:Jump}
\end{figure}

\end{example}

\begin{example}(continued from Example \ref{example:diamond_graph}) \label{example:diamond_graph_contd} Consider the coupled wave equations on the diamond network described from Example \ref{example:diamond_graph} with different velocities of propagation  $c_{e_i}$ on each edge $e_i \in E$ with $i=1,...,7$. 

To illustrate the construction of the randomized approximation, introduce four subsets as
$E_1 = \{1, 2, 3\}, E_2 = \{2, 4, 5\}, E_3 = \{ 3,4,6\}, E_4 = \{5,6,7\}$, which satisfy $\cup_{\omega=1}^4 S_\omega = \{1,...,7\}$. Each subset $E_\omega$ generates a tripod (star-shaped) subgraph as shown in Figure \ref{fig:ram_diamond}. On each time interval, one of these four subgraphs is chosen with the same possibility $p_\omega = \frac{1}{4} (\omega = 1,...,4)$. The solution on the tripod network is easier to compute than the solution on the original network because the subnetworks do not contain loops. 

\begin{figure}[h]
    \centering
    \includegraphics[trim={0 0 6.3cm 0},clip,width=0.7\textwidth]{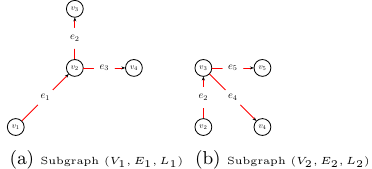}
    
    \includegraphics[trim={6.3cm 0 0 0},clip,width=0.7\textwidth]{Fig3.pdf}
    \caption{A decomposition of the diamond network in Figure \ref{fig:diamond-graph} into four (overlapping) subgraphs}
    \label{fig:ram_diamond}
\end{figure}

In this example, each edge has the different probability of being selected. To be specific, $\pi_1 = \pi_7 = \frac{1}{4}$ ($e_1$ and $e_7$ only appear in one of the subnetworks), while $\pi_{e_2} = \pi_{e_3} = \pi_{e_4} = \pi_{e_5} =\pi_{e_6} = \frac{1}{2}$ (these edges appear in two of the subnetworks). Thus, for $t\in (t_{k-1}, t_k]$, the velocities of propagation in the randomized network are enlarged differently on each edge $e_i$ based on these probabilities. In particular, 
\begin{equation}
    c_{h,e_i} = \left\{\begin{aligned}
        &4 c_{e_i}, \quad & \text{for} \quad e_i\in \{e_1,e_7\} \cap  E_{\omega_k}\\
        &2 c_{e_i}, \quad &\text{for} \quad e_i\in \{ e_2,e_3,e_4,e_5,e_6\} \cap  E_{\omega_k}\\
        &0, \quad  &\text{for} e_i\notin E_{\omega_k}.
    \end{aligned}\right.
\end{equation}
Based on this randomized velocity field, the corresponding randomized dynamical system \eqref{eq:RD} and its optimal control can be computed and will approximate the solution of the original system \eqref{eq:D} when the length of the time intervals $h$ is sufficiently small. 
\end{example} 

\subsection{Main results}
This section summarizes the two main results of this paper. The proofs are postponed until Sections \ref{sec:dynamics} and \ref{sec:control}. 

The first main result shows that the solutions $w_{h+}^{e_i}(\bomega,t,x)$, $w_{h-}^{e_i}(\bomega,t,x)$, and $y_h^{e_i}(\bomega,t,x)$ to the randomized dynamics \eqref{eq:RD} approach the solutions $w_{+}^{e_i}(t,x)$, $w_{-}^{e_i}(t,x)$, and $y^{e_i}(t,x)$ to the original dynamics \eqref{eq:D}, respectively, for $h \rightarrow 0$. 

\begin{theorem} \label{thm:dyn2}
    Let $T > 0$ and the initial conditions $(\mathbf{y}_0, \mathbf{y}_1)$ and the boundary control $\mathbf{u}$ be such that 1) they satisfy the BCs in \eqref{eq:D}, 2) $y_{0,x}^{e_i}$ and $y_1^{e_i}$ are Lipschitz in $x$ for every edge $e_i \in E$, and 3) $\mathbf{u} $ is Lipschitz in $t$. Then there exists a constant $C$ that does not depend on $h$ such that for any $e_i \in E$, $x \in (0, \ell_{e_i})$ and $0 \leq t \leq T$
$$
        \mathbb{E}[|w_{h-}^{e_i}(t,x) - w_-^{e_i}(t,x)|^2 + |w_{h+}^{e_i}(t,x) - w_+^{e_i}(t,x)|^2 + |y_h^{e_i}(t,x) - y^{e_i}(t,x)|^2] \leq C h. 
 $$
\end{theorem}
The proof of Theorem \ref{thm:dyn2} will be given in Section \ref{sec:dynamics}. 

\begin{remark}
The convergence in expectation in Theorem \ref{thm:dyn2} implies the convergence in probability. This follows from Markov's inequality, which states that for any random variable $X\geq 0$ and $a > 0$
    \begin{equation}
        \mathbb{P}[X \geq a] \leq \frac{\mathbb{E}[X]}{a}. 
    \end{equation}
This also means that
$$
        \mathbb{P}[|w_{h-}^{e_i}(t,x) - w_-^{e_i}(t,x)|^2 + |w_{h+}^{e_i}(t,x) - w_+^{e_i}(t,x)|^2 + |y_h^{e_i}(t,x) - y^{e_i}(t,x)|^2 \geq \varepsilon] \leq \frac{C h}{\varepsilon}. 
$$
In other words, the probability that the distance between the randomized solution of \eqref{eq:RD} and the solution to the original problem \eqref{eq:D} is more than any $\varepsilon > 0$ can be made arbitrarily small by choosing the spacing of the temporal grid $h$ sufficiently small. 
\end{remark}

For the randomized optimal control problem we furthermore have the following convergence result.

\begin{theorem} \label{thm:contr}
When the initial conditions $(\mathbf{y}_0, \mathbf{y}_1)$ are as in Theorem \ref{thm:dyn2}, there exists a constant $C$ that does not depend on $h$ such that
\begin{equation}
    \lim_{h \rightarrow 0} \mathbb{E}[|\mathbf{u}^*_h - \mathbf{u}^*|_{H^2(0,T; \mathbb{R}^{|V_C|})}^2] \leq Ch. 
\end{equation}
    
\end{theorem}
The proof of Theorem \ref{thm:contr} will be given in Section \ref{sec:control}. 

Note that Theorem \ref{thm:contr} shows that the probability that the distance between the optimal control for the randomized system \eqref{eq:RD} and the optimal control for the original system \eqref{eq:D} is larger than any $\varepsilon > 0$ can be made arbitrarily small by choosing $h$ sufficiently small.

\section{Preliminaries} \label{sec:preliminaries}

\subsection{Preliminaries for \texorpdfstring{\eqref{eq:D}}{(D)} \label{ssec:preliminariesD}}


Just as before, write $\Omega = \Pi_{e_i\in E} (0, l_{e_i})$ and consider
\begin{align*}
&\mathcal{H} = L^2(\Omega), &&\| \mathbf{y} \|_{\mathcal{H}} = \sum_{e_i \in E} |y^{e_i}|_{L^2(0,\ell_{e_i})}, \\
&\mathcal{H}^1 = \{ \mathbf{y} \in H^1(\Omega) \mid y^{e_i}(v_j) = y^{e_k}(v_j), e_i, e_k \in E(v_j) \}, && \| \mathbf{y} \|_{\mathcal{H}^1} = \sum_{e_i \in E} |y^{e_i}|_{H^1(0,\ell_{e_i})}. 
\end{align*}

For the forward dynamics \eqref{eq.system} there is the following classical existence result , see e.g. \cite[Theorem 2.53]{coron2007}.
\begin{theorem}[Wellposedness]
\label{thm. wellposedness} 
For $T>0$, initial data $(\mathbf{y}_0, \mathbf{y}_1) \in  \mathcal{H}^1 \times \mathcal{H}$ and boundary controls $\mathbf{u} \in L^2(0,T; \mathbb{R}^{|V_C|})$, there exists a unique solution $\mathbf{y} \in C(0,T; \mathcal{H}^1) \bigcap C^1 (0,T; \mathcal{H})$ to the initial boundary value problem \eqref{eq:D}. This solution satisfies the estimate
\begin{equation}
\| \mathbf y \|_{C(0,T; \mathcal{H}^1)} + \| \mathbf y_t \|_{C(0,T; \mathcal{H})}\le C_{\mathrm{D}}( \|\mathbf y_0\|_{\mathcal{H}^1} + \|\mathbf y_1\|_{\mathcal{H}} + | \mathbf{u} |_{L^2(0,T; \mathbb{R}^{|V_C|})}), 
\label{eq:wellposedness1}
\end{equation}
where the constant $C_{\mathrm{D}}$ is independent of $\mathbf y_0, \mathbf y_1$ and $\mathbf u$. 

\end{theorem}


The following proofs will use that the Riemann invariants $w_+^{e_i}(t,x)$ and $w_-^{e_i}(t,x)$ are Lipschitz in their arguments $t$ and $x$. The following notation is used. For Lipschitz-continuous functions $y_0^{e_i}(x)$ or $\mathbf{u}(t)$ defined on an edge $e_i \in E$ or depending on time, let $\mathrm{Lip}(y_0^{e_i})$ and $\mathrm{Lip}(\mathbf{u})$ denote their Lipschitz constants, respectively. For a function $\mathbf{y}_0$ on $\Omega$, define
\begin{equation}
\mathrm{Lip}(\mathbf{y}_0) = \max_{e_i \in E} \mathrm{Lip}(y_0^{e_i}). 
\end{equation}
The following Theorem can be obtained easily by tracing back the solutions along characteristics. 
\begin{theorem}[Lipschitz solutions] \label{thm:Lipschitz}
Let $T > 0$. If the initial conditions $(\mathbf{y}_0, \mathbf{y}_1)$ and the boundary control $\mathbf{u}$ 1) satisfy the boundary conditions in \eqref{eq:D}, 2) $y_{0,x}^{e_i}$ and $y_1^{e_i}$ are Lipschitz in $x$ for all $e_i \in E$ , and 3) $\mathbf{u}$ is Lipschitz in $t$, then $w_{+}^{e_i}(t,x)$ and $w_{-}^{e_i}(t,x)$ are Lipschitz in their arguments $t$ and $x$. In particular, there exists a constant $C_{\mathrm{Lip}}$, depending on $T$, such that for all $e_i \in E$, $x'', x' \in (0, \ell_{e_i})$, and $t'' ,t' \in [0,T]$
\begin{align}
|w^{e_i}_\pm(t'', & x'') - w^{e_i}_\pm(t', x')| \nonumber \\ 
& \qquad \qquad \leq C_{\mathrm{Lip}} (\mathrm{Lip}(\mathbf{y}_{0,x}) + \mathrm{Lip}(\mathbf{y}_{1}) + \mathrm{Lip}(\mathbf{u})) \left(  |x'' - x'| + \bar{c} |t'' - t'| \right),  \label{eq:thm_Lipschitz}
\end{align}
where $\bar{c} = \max_{e_i \in E} c_{e_i}$. 
\end{theorem}

Because elements of the proof of Theorem \ref{thm:Lipschitz} also appear in the convergence proof of our proposed algorithm, the proof of Theorem \ref{thm:Lipschitz} is given in Appendix \ref{app:thm_Lipschitz}.

\begin{theorem}[Existence and uniqueness of optimal controls] \label{thm:OCP} 
The optimal control problem \eqref{eq:OCP} has a unique solution $\mathbf{u}^* \in U_{ad}$ for any $(\mathbf{y}_0, \mathbf{y}_1) \in \mathcal{H}^1 \times \mathcal{H}$.  
\end{theorem}

The theorem follows almost directly from the strict convexity of the cost functional, see e.g.\ \cite{minoux1986}. Note that any control $\mathbf{u} \in H^2(0,T; \mathbb{R}^{|V_C|})$ is Lipschitz. 

\subsection{Probabilistic Framework}

The probability space for the considered randomized dynamical system is
\begin{equation}
    \mathcal{P} = \{1,2, \ldots, 2^{|E|} \}^K. 
\end{equation}
The chosen probabilities $p_\omega$ with $\omega \in \{1,2, \ldots, 2^{|E|}\}$ assign a probability to each point $\bomega = (\omega_1, \omega_2, \ldots, \omega_K)$ in $\mathcal{P}$ according to
\begin{equation}
    p(\bomega) = \prod_{k=1}^K p_{\omega_k}. 
\end{equation}
The expected value of a random variable $X(\bomega)$ is defined as
\begin{equation}
    \mathbb{E}[X] = \sum_{\bomega \in \mathcal{P}} X(\bomega) p(\bomega). 
\end{equation}
Note that
\begin{equation}
    \mathbb{E}[X] = \sum_{\omega_1=1}^{2^{|E|}} 
    \sum_{\omega_2=1}^{2^{|E|}}
    \cdots
    \sum_{\omega_K=1}^{2^{|E|}}
    X(\omega_1, \omega_2, \ldots, \omega_K) p_{\omega_1}p_{\omega_2} \cdots p_{\omega_K}. 
\end{equation}
Therefore, \eqref{eq:sump} shows that if $X = X(\omega_k)$ only depends on one element $\omega_k$ of $\bomega$,
\begin{equation}
    \mathbb{E}[X] =  
    \sum_{\omega_2=k}^{2^{|E|}}
    X(\omega_k) p_{\omega_k}, 
\end{equation}
where it has been used that the $p_\omega$'s satisfy \eqref{eq:sump}. 

    It is then easy to see that (with $t \in (t_{k-1}, t_k)$ so that $c_{h,e_i}(\bomega,t) = c_{h,e_i}(\omega_k,t)$)
    \begin{equation}
        \mathbb{E}[c_{h,e_i}(t)] := \sum_{\omega_k = 1}^{2^{|E|}} c_{h,e_i}(\omega_k, t) p_{\omega_k} = \frac{c_{e_i}}{\pi_{e_i}} \sum_{\omega
        _k \in \{ \omega \mid e_i \in S_{\omega} \}} p_{\omega_k} = c_{e_i}, 
    \end{equation}
    where the last identity follows from \eqref{eq:def_pi}. To measure the deviation from $c_{h,e_i}(\bomega,t)$ from $c_{e_i}$, consider (for $t \in [t_{k-1}, t_k]$)
    \begin{align}
        \mathbb{E}[|c_{h,e_i}(t) - c_{e_i}|^2] &= \sum_{\omega_k = 1}^{2^M} (c_{h,e_i}(\omega_k,t) - c_{e_i})^2 p_{\omega_k} \nonumber \\ 
        &= c_{e_i}^2 \sum_{\omega=1}^{2^M} \left(\frac{\chi_{e_i}(\omega)}{\pi_{e_i}} - 1 \right)^2 =: \mathrm{Var}[c_{h,e_i}]. \label{eq:Varc}
    \end{align}
    with $\chi_{e_i}$ as in \eqref{eq:def_pi}. Also observe that the third expression is in fact independent of the considered time instant $t$ and is denoted by $\mathrm{Var}[c_{h,e_i}]$.

\subsection{Preliminaries for \texorpdfstring{\eqref{eq:RD}}{(RD)} \label{ssec:preliminariesRD}}
For the randomized system \eqref{eq:RD}, there is the following existence and uniqueness result, which is an analogue to Theorem \ref{thm. wellposedness} in Subsection \ref{ssec:preliminariesD} for the original system \eqref{eq:D} but with lower regularity because of the discontinuities introduced by the randomization see e.g. Example \ref{example:jump}. 

\begin{theorem}[Wellposedness for random networks]
\label{thm. wellposednessh} 
Let $T>0$ and consider initial data $(\mathbf{y}_0, \mathbf{y}_1) \in \mathcal{H}^1 \times \mathcal{H}$. For every $\bomega \in \mathcal{P}$ and boundary control $\mathbf{u} \in H^2(0,T; \mathbb{R}^{|V_C|})$, there exists a unique solution $\mathbf{y}_h(\bomega)$ to the initial boundary value problem \eqref{eq:RD}. This solution satisfies the estimate
$$
\| \mathbf y_h(\bomega) \|_{C(0,T; \mathcal{H})} \le C_{\mathrm{RD}}( \|\mathbf y_0\|_{\mathcal{H}^1} + \|\mathbf y_1\|_{\mathcal{H}} + | \mathbf{u} |_{L^2(0,T;\mathbb{R}^{|V_C|})}),
$$
where the constant $C_{\mathrm{RD}}$ is independent of $\mathbf y_0, \mathbf y_1$, $\mathbf u$, $\bomega$, and the time grid $t_0, t_1, \ldots, t_K$. 
\end{theorem}

Furthermore, there is the following existence and uniqueness result for the optimal controls of \eqref{eq:OCPh}. 

\begin{theorem}[Randomized Optimal Control] \label{thm:OCPh}
Let $T > 0$ and $(\mathbf{y}_0, \mathbf{y}_1) \in \mathcal{H}^1 \times \mathcal{H}$. For every $\bomega \in \{1,2, \ldots, 2^{|E|} \}^K$, the optimal control problem \eqref{eq:OCPh} has a unique solution $\mathbf{u}_h^*(\bomega) \in U_{ad}$. Furthermore, 
there exist constants $C_1$ and $C_2$ independent of $\bomega$ and the time grid $t_0, t_1, \ldots, t_K$ such that 
\begin{equation}
|\mathbf{u}^*_h(\bomega)|_{H^2(0,T; \mathbb{R}^{|V_C|})} \leq C_1, \qquad \qquad 
\mathrm{Lip}(\mathbf{u}^*_h(\bomega)) \leq C_2. 
\end{equation}
\end{theorem}

\begin{proof}
The existence and uniqueness of the optimal control follow just as in Theorem \ref{thm:OCP} from the strict convexity of the cost functional, see e.g.\ \cite{minoux1986}. 

The uniform bound on $|\mathbf{u}^*_h(\bomega)|_{H^2(0,T; \mathbb{R}^{|V_C|})}$ follows by noting that
\begin{align*}
\frac{\alpha}{2} |\mathbf{u}^*_h(\bomega)&|_{H^2(0,T; \mathbb{R}^{|V_C|})}^2 \leq J_h(\bomega, \mathbf{u}^*_h(\bomega)) \leq J_h(\bomega, \mathbf{0}) \\
&= \frac{1}{2} \| \mathbf{y}_h(\bomega) - \mathbf{y}_d \|^2_{L^2(Q)} \leq \frac{1}{2}\left( TC_{\mathrm{RD}}(\| \mathbf{y}_0 \|_{\mathcal{H}^1} + \| \mathbf{y}_1 \|_{\mathcal{H}}) + \| \mathrm{y}_d \|_{L^2} \right)^2, 
\end{align*}
where the first inequality uses the definition of $J_h(\bomega,\mathbf{u})$ in \eqref{eq:OCPh} and the second inequality follows because $\mathbf{u}^*_h(\bomega)$ is the minimizer of $J_h(\bomega,\cdot)$. On the second line, $\mathbf{y}_h(\bomega)$ denotes the solution to \eqref{eq:RD} and \eqref{eq:yhintwh} resulting from zero control, and the last inequality uses the triangle inequality and Theorem \ref{thm. wellposednessh}. 

The uniform bound on $\mathrm{Lip}(\mathbf{u}^*_h(\bomega))$ directly from the Sobolev inequalities which show that $H^2(0,T; \mathbb{R}^{|V_C|})$ is compactly embedded in $W^{1,\infty}(0,T; \mathbb{R}^{|V_C|})$, see e.g.\ \cite{evans1998}. 
\end{proof}

\section{Convergence of the Forward Dynamics \label{sec:dynamics}}
\subsection{Convergence of the Characteristics}

Characteristics are important to describe the wave propagation along networks and play a central role in the proof of Theorem \ref{thm:dyn2}. Let $\xi^{e_i}_{\pm}(s;t,x)$ and $\xi^{e_i}_{h\pm}(\bomega,s;t,x)$ denote the characteristic curves passing through $(t, x)$ which satisfy by definition
\begin{multline}
  \frac{\mathrm d \xi^{e_i}_{\pm}}{\mathrm ds}(s;t,x) = \pm c_{e_i}, \qquad \frac{\mathrm d \xi^{e_i}_{h\pm}}{\mathrm ds}(\bomega,s;t,x) = \pm c_{h,e_i}(\bomega,s), \\ \xi^{e_i}_\pm(t;t,x) = \xi^{e_i}_{h\pm}(\bomega,t,t,x)= x. 
\end{multline}
It is easy to see that $\xi^{e_i}_{\pm}(s;t,x)$ and $\xi^{e_i}_{h\pm}(\bomega,s;t,x)$ are given by
\begin{equation}
    \xi^{e_i}_{\pm}(s;t,x) = x \mp c_{e_i} (t-s), \qquad \qquad 
    \xi^{e_i}_{h\pm}(\bomega,s;t,x) = x \mp \int_s^t c_{h,e_i}(\bomega,s') \ \mathrm{d}s'. \label{eq:xi_sol}
\end{equation}
We consider the functions in \eqref{eq:xi_sol} for all $0 \leq s \leq t$ and $x \in (0, \ell_{e_i})$. In particular, in the following the values of $\xi^{e_i}_{\pm}(s;t,x)$ and $\xi^{e_i}_{h\pm}(\bomega,s;t,x)$ may not be inside $(0,\ell_{e_i})$. 


The following lemma shows that the randomized characteristics converge to the characteristics of the original system when $h$ approaches zero. 
\begin{lemma} \label{lem:conv_xi}
Let $\mathrm{Var}[c_{h,e_i}]$ be as in \eqref{eq:Varc}. Then for all $0 \leq s \leq t$,
    \begin{equation}
        \mathbb{E}\left[|\xi^{e_i}_{h\pm}(s; t) - \xi^{e_i}_\pm(s; t)|^2_{L^\infty(0,\ell_{e_i})} \right] \leq h(t-s)\mathrm{Var}[c_{h,e_i}]. \label{eq:lem_xi1}
    \end{equation}
Furthermore, there exists a constant $C_1$ such that
 \begin{equation}
        \mathbb{E}\left[|\xi^{e_i}_{h\pm}(s; t) - \xi^{e_i}_\pm(s; t)|^4_{L^\infty(0,\ell_{e_i})} \right] \leq C_1 h^2(t-s)^2. \label{eq:lem_xi2}
    \end{equation}
\end{lemma}

\begin{proof} Equation \eqref{eq:xi_sol} shows that
    \begin{equation}
        \xi^{e_i}_{h\pm}(\bomega,s;t,x) - \xi^{e_i}_{\pm}(s;t,x) = \mp \int_s^t \left( c_{h,e_i}(\bomega,s') - c_{e_i} \right) \ \mathrm{d}s'. 
    \end{equation}
    Let $t' = t'_0 < t'_1 < \ldots < t'_{L-1} < t'_L = t$ be the restriction of the timegrid $\{ t_k \}_k$ to the time interval $[s, t]$ and note that
    \begin{align}
        &\left| \xi^{e_i}_{h\pm}(\bomega,s;t) - \xi^{e_i}_{\pm}(s;t) \right|^2_{L^\infty(0,\ell_{e_i})} = \left( \int_{s}^{t} \left( c_{h,e_i}(\bomega,s') - c_{e_i} \right) \ \mathrm{d}s' \right)^2 \nonumber \\
        &\qquad = \sum_{k=1}^L\sum_{\ell=1}^L \left( \int_{t'_{k-1}}^{t'_k} \left( c_{h,e_i}(\bomega,s') - c_{e_i} \right) \ \mathrm{d}s' \right) \left( \int_{t'_{\ell-1}}^{t'_\ell} \left( c_{h,e_i}(\bomega,s') - c_{e_i} \right) \ \mathrm{d}s' \right). 
    \end{align}
    When taking the expectation, the terms with $k \neq \ell$ vanish because $\mathbb{E}[c_{h,e_i}(s)] = c_{e_i}$ and $c_{h,e_i}(\bomega,s')$ and $c_{h,e_i}(\bomega,s'')$ are independent when $s'$ and $s''$ are in different sub-intervals $[t'_{k-1}, t'_k]$ and $[t'_{k'-1}, t'_{k'}]$, respectively. Therefore, 
    \begin{align}
        &\mathbb{E}\left[ \left|\xi^{e_i}_{h\pm}(s;t) - \xi^{e_i}_{\pm}(s;t) \right|^2_{L^\infty(0,\ell_{e_i})} \right] =  \sum_{\ell=1}^L \mathbb{E}\left[ \left( \int_{t'_{\ell-1}}^{t'_\ell} \left( c_{h,e_i}(s') - c_{e_i} \right) \ \mathrm{d}s' \right)^2 \right] \nonumber \\
        &\leq \sum_{\ell=1}^L \mathbb{E}\left[ (t'_\ell-t'_{\ell-1}) \int_{t'_{\ell-1}}^{t'_\ell} \left( c_{h,e_i}(s') - c_{e_i} \right)^2 \ \mathrm{d}s' \right] \leq h \sum_{\ell=1}^L \int_{t'_{\ell-1}}^{t'_\ell} \mathrm{Var}[c_{h,e_i}] \ \mathrm{d}s \nonumber \\ 
        &= h(t-s) \mathrm{Var}[c_{h,e_i}], \label{eq:lem_xi_step3}
    \end{align}
    where the inequality follows from Cauchy-Schwarz and the second to last equality uses \eqref{eq:Varc}. This completes the proof of \eqref{eq:lem_xi1}. 

For \eqref{eq:lem_xi2}, write 
\begin{equation}
    F_{k}^{e_i}(\bomega) = \int_{t'_{k-1}}^{t'_k} (c_{h,e_i}(\bomega,s') - c_{e_i}) \ \mathrm{d}s' = (t'_k - t'_{k-1}) (c_{h,e_i}(\omega_k,s) - c_{e_i}). \label{eq:lem_defFk}
\end{equation}
and observe that
\begin{align}
\mathbb{E}\left[ \left|\xi^{e_i}_{h\pm}(s;t) - \xi^{e_i}_{\pm}(s;t) \right|^4_{L^\infty(0,\ell_{e_i})} \right] &= \mathbb{E}\left[ \left( \int_{s}^{t} (c_{h,e_i}(s') - c_{e_i}) \ \mathrm{d}s' \right)^4 \right] \nonumber \\ 
&= \sum_{k_1,k_2,k_3,k_4=1}^{L} \mathbb{E}[ F_{k_1}^{e_i} F_{k_2}^{e_i} F_{k_3}^{e_i} F_{k_4}^{e_i} ]
\end{align}
Because $F^{e_i}_{k}$ and $F^{e_i}_{\ell}$ are independent when $k \neq \ell$ and $\mathbb{E}[F^{e_i}_k] = 0$, the expected value $\mathbb{E}[F^{e_i}_{k_1}F^{e_i}_{k_2}F^{e_i}_{k_3}F^{e_i}_{k_4}]$ vanishes unless there is a $n \in \{ 2,3,4 \}$ such that $k_n = k_1$. Because the indices $k_1$, $k_2$, $k_3$, and $k_4$ are interchangeable, 
\begin{equation}
    \mathbb{E}\left[ \left|\xi^{e_i}_{h\pm}(t';t) - \xi^{e_i}_{\pm}(t';t) \right|^4_{L^\infty(0,\ell_{e_i})} \right]
\leq 3 \sum_{k_1,k_3,k_4=1}^L \mathbb{E}[(F^{e_i}_{k_1})^2 F^{e_i}_{k_3} F^{e_i}_{k_4}].  
\end{equation}
Because $c_{h,e_i}(\omega_k,s) - c_{e_i}$ is either $c_{e_i}$ or $c_{e_i} / \pi_{e_i} - c_{e_i}$, looking back at \eqref{eq:lem_defFk} shows that
\begin{equation}
(F_k^{e_i}(\bomega))^2 \leq C_0^2 (t'_{k} - t'_{k-1})^2, \qquad \qquad C_0 = \max_{e_i \in E} \left(  c_{e_i} \max \left\{1, \left| \frac{\pi_{e_i} - 1}{\pi_{e_i}} \right|  \right\} \right). 
\end{equation}
Therefore, 
\begin{align}
    &\mathbb{E}\left[ \left|\xi^{e_i}_{h\pm}(s;t) - \xi^{e_i}_{\pm}(s;t) \right|^4_{L^\infty(0,\ell_{e_i})} \right] 
\leq 3 \sum_{k_1,k_3,k_4=1}^L C_0^2 (t'_{k_1} - t'_{k_1-1})^2 \mathbb{E}[F^{e_i}_{k_3} F^{e_i}_{k_4}] \nonumber \\
&\qquad \leq 3 C_0^2 h (t-s) \mathbb{E}\left[\left( \int_{t}^{t'} (c_{h,e_i}(s) - c_{e_i}) \ \mathrm{d}s \right)^2\right] \leq 3 C_0^2 h^2 (t-s)^2 \mathrm{Var}[c_{e_i}], 
\end{align}
where the last inequality follows from \eqref{eq:lem_xi_step3}. This is \eqref{eq:lem_xi2} with $C_1 = 3C_0^2\mathrm{Var}[c_{h,e_i}]$. 
\end{proof}

The time instants $s \in [0,t]$ at which the characteristic curves $s \mapsto (s, \xi^{e_i}_\pm(s; t,x))$ leave the domain $(0,t) \times (0,\ell_{e_i})$ are denoted by $t_{\pm,\mathrm{in}}^{e_i}(t,x)$ and the time instants at which the characteristic curves $s \mapsto (s, \xi^{e_i}_{h\pm}(\bomega,s;t,x))$ leave $(0,t) \times (0,\ell_{e_i})$ are denoted by $t_{h\pm,\mathrm{in}}^{e_i}(\bomega,t,x)$. 

For some $t' \in [0,t)$, it is easy so check that $\max\{ t_{\pm,\mathrm{in}}^{e_i}(t), t' \}$ and $\max\{ t_{h\pm,\mathrm{in}}^{e_i}(t), t' \}$ are the time instants $s$ at which the characteristic curves $s \mapsto (s, \xi^{e_i}_\pm(s; t,x))$ and $s \mapsto (s, \xi^{e_i}_{h\pm}(\bomega,s;t,x))$ leave the domain $(t',t) \times (0,\ell_{e_i})$, respectively. The time instant $t'$ will appear naturally in the proof of Theorem \ref{thm:dyn2}. 

The following lemma shows that $\max\{ t_{h\pm,\mathrm{in}}^{e_i}(t), t' \}$ converges in expectation to $\max\{ t_{\pm,\mathrm{in}}^{e_i}(t), t' \}$ when $h$ is sufficiently small. 

\begin{lemma} \label{lem:conv_tini}
There exists a constant $C_2$ independent of $h$ s.t.\ for all $0 \leq t' \leq t \leq T$
\begin{equation}
\mathbb{E}[|\max\{ t_{h\pm,\mathrm{in}}^{e_i}(t), t' \} - \max\{ t_{\pm,\mathrm{in}}^{e_i}(t), t' \} |^2_{L^\infty(0,\ell_{e_i})}] \leq C_2 h (t - t'). 
\label{eq:lem_conv_tini}
\end{equation}
\end{lemma}

\begin{proof}
We only prove the result for $t_{h+,\mathrm{in}}^{e_i}(\bomega,t,x)$ and $t_{+,\mathrm{in}}^{e_i}(t,x)$. The result for $t_{h-,\mathrm{in}}^{e_i}(\bomega,t,x)$ and $t_{h-,\mathrm{in}}^{e_i}(\bomega,t,x)$ can be proven in a similar way. 
For brevity, we write $t^*_h(\bomega,t,x)$ and $t^*(t,x)$ for $\max\{ t_{h+,\mathrm{in}}^{e_i}(\bomega,t,x), t' \}$ and $\max\{ t_{+,\mathrm{in}}^{e_i}(t,x), t' \}$, respectively. 
By definition of the Lebesque integral
\begin{equation}
\mathbb{E}[|t^*_h(t) -t^*(t)|^2_{L^\infty(0,\ell_{e_i})}] = \int_0^\infty F(\tau) \ \mathrm{d}\tau,  \quad F(\tau) = \mathbb{P}[|t^*_h(t) -t^*(t)|^2_{L^\infty(0,\ell_{e_i})} > \tau] . \label{eq:exp_cumF}
\end{equation}
An upperbound for $F(\tau)$ thus leads to an upperbound for $\mathbb{E}[|t^*_h(t) -t^*(t)|^2_{L^\infty(0,\ell_{e_i})}]$.

When $|t^*_h(\bomega,t) - t^*(t)|_{L^\infty(0,\ell_{e_i})} > \sqrt{\tau}$, there exists an $x \in [0,\ell_{e_i}]$ such that either $t^*_h(\bomega,t,x) > t^*(t,x) + \sqrt{\tau}$ or  $t^*_h(\bomega,t,x) < t^*(t,x) - \sqrt{\tau}$. Both cases are depicted in Figure \ref{fig:tini_conv}. 

Figure \ref{fig:tini_conv1} shows the case where $t^*_h(\bomega,t,x) > t^*(t,x) + \sqrt{\tau}$. Because $t^*_h(\bomega,t,x) \in [t',t]$, this can only occur when $s_+ := t^*(t,x) + \sqrt{\tau} \in [t',t]$. Furthermore, the characteristic $\xi_{h+}^{e_i}(\bomega,s; t,x) < 0$ for all $s < t^*_h(\bomega,t,x)$, so $\xi^{e_i}_{h+}(\bomega, s_+; t,x) < 0$ in this case. On the other hand, because $\xi^{e_i}_+(t^*(t,x); t,x) \geq 0$ and the characteristic $\xi_+^{e_i}(\bomega,s;t,x)$ travels with velocity $c_{e_i}$, it follows that $\xi^{e_i}_+(t^*(t,x) + \sqrt{\tau}; t,x) \geq c_{e_i} \sqrt{\tau}$. Therefore,
\begin{equation}
t^*_h(\bomega,t,x) > t^*(t,x) + \sqrt{\tau} \quad \Rightarrow \quad
|\xi^{e_i}_{h+}(\bomega,s_+; t,x) - \xi^{e_i}_+(s_+; t,x)| > c_{e_i} \sqrt{\tau}. \label{eq:lem_conv_tini_step2}
\end{equation}
Next, consider Figure \ref{fig:tini_conv2} which depicts the case where $t^*_h(\bomega,t,x) < t^*(t,x) - \sqrt{\tau}$. Because $t^*_h(\bomega,t,x) \in [t',t]$, this can only occur when $s_- := t^*(t,x) - \sqrt{\tau} \in [t',t]$. It can be checked similarly as before that $\xi_{h+}^{e_i}(\bomega,s_-; t,x) > 0$ and $\xi^{e_i}_h(s_-; t,x) \leq -c_{e_i} \sqrt{\tau}$ in this case. Therefore,
\begin{equation}
t^*_h(\bomega,t,x) < t^*(t,x) - \sqrt{\tau} \quad \Rightarrow \quad 
|\xi_{h+}^{e_i}(\bomega,s_-; t,x) - \xi^{e_i}_+(s_-; t,x)| > c_{e_i}\sqrt{\tau}. \label{eq:lem_conv_tini_step3}
\end{equation}

\begin{figure}
\centering

\subfloat[$t_{h+,\mathrm{in}}^{e_i}(\bomega;t,x) > t^{e_i}_{+,\mathrm{in}}(t,x)$ \label{fig:tini_conv1}]{
\includegraphics[width=0.45\textwidth]{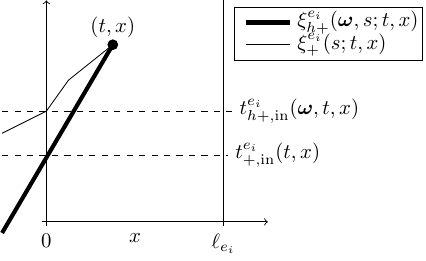}
}
~
\subfloat[$t_{h+,\mathrm{in}}^{e_i}(\bomega,t,x) < t^{e_i}_{+,\mathrm{in}}(t,x)$ \label{fig:tini_conv2}]{
\includegraphics[width=0.45\textwidth]{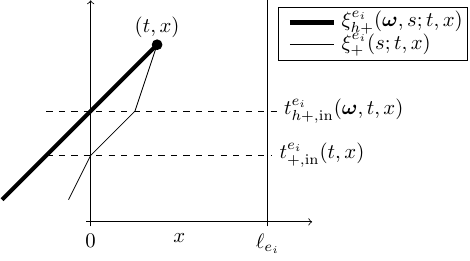}
}
\caption{The characteristics $\xi^{e_i}_+(s;t,x)$ and $\xi_{h+}^{e_i}(\bomega,s;t,x)$ with a positive characteristic velocity cross the boundary $x = 0$ at times $t^{e_i}_{+,\mathrm{in}}(t,x)$ and $t^{e_i}_{h+,\mathrm{in}}(\bomega,t,x)$, respectively. }
\label{fig:tini_conv}
\end{figure}

Combining \eqref{eq:lem_conv_tini_step2} and \eqref{eq:lem_conv_tini_step3} and raising the inequalities on the LHS to the power two and raising the inequalities on the RHS to the power four shows that
\begin{align}
|t^*_h(\bomega,t) - t^*(t)|_{L^\infty(0,\ell_{e_i})}^2 > \tau & \quad \Rightarrow |\xi^{e_i}_{h+}(\bomega,s_+; t) - \xi^{e_i}_+(s_+; t)|_{L^\infty(0,\ell_{e_i})}^4 > c_{e_i}^4 \tau^2 \nonumber \\
\mathrm{or}\ 
& |\xi^{e_i}_{h+}(\bomega,s_-; t) - \xi^{e_i}_+(s_-; t)|_{L^\infty(0,\ell_{e_i})}^4 > c_{e_i}^4 \tau^2. 
\end{align}

By Markov's inequality and \eqref{eq:lem_xi2} from Lemma \ref{lem:conv_xi}, 
\begin{align}
    \mathbb{P}[|t^*_h(t) & -t^*(t)|^2_{L^\infty(0,\ell_{e_i})} > \tau] \leq \mathbb{P}[|\xi^{e_i}_{h+}(\bomega,s_+; t) - \xi^{e_i}_+(s_+; t)|_{L^\infty(0,\ell_{e_i})}^4 > c_{e_i}^4 \tau^2] \nonumber \\
    &\qquad\qquad\qquad\qquad\qquad\qquad +\mathbb{P}[|\xi^{e_i}_{h+}(\bomega,s_-; t) - \xi^{e_i}_+(s_-; t)|_{L^\infty(0,\ell_{e_i})}^4 > c_{e_i}^4 \tau^2] \nonumber \\
    &\leq \frac{\mathbb{E}[|\xi^{e_i}_{h+}(s_+; t) - \xi^{e_i}_+(s_+; t)|^4] + \mathbb{E}[|\xi^{e_i}_{h+}(s_-; t) - \xi^{e_i}_+(s_-; t)|^4]}{c_{e_i}^4\tau^2}  \nonumber \\
    &\leq \frac{C_1h^2(t-s_+)^2 + C_1h^2(t-s_-)^2}{c_{e_i}^4\tau^2} \leq \frac{2C_1h^2(t-t')^2}{c_{e_i}^4 \tau^2},
\end{align}
where it has been used that $s_\pm \in [t', t]$ in . Because a probability can never exceed one, it follows that
\begin{equation}
F(\tau) \leq \min\left\{1, \frac{G^2h^2}{\tau^2} \right\}, \qquad \qquad G := \frac{(t-t')\sqrt{2C_1}}{c_i^2}.
\end{equation}
Inserting this estimate into \eqref{eq:exp_cumF}, it follows that
\begin{equation}
\mathbb{E}[|t^*_h - t^*|^2] \leq \int_0^\infty \min\left\{1, \frac{G^2 h^2}{\tau^2} \right\} \ \mathrm{d}\tau = \int_0^{G h} \ \mathrm{d}\tau + \int_{G h}^\infty \frac{G^2 h^2}{\tau^2} \ \mathrm{d}\tau = \frac{3G h}{2} . 
\end{equation}
This yields \eqref{eq:lem_conv_tini} with $C_2 = \tfrac{3}{2}\sqrt{2C_1} / \underline{c}^2$, where $\underline{c} = \min_i c_{e_i} > 0$. 
\end{proof}

\subsection{Proof of Theorem \ref{thm:dyn2}}

The following lemma establishes the convergence in the Riemann variables in Theorem \ref{thm:dyn2}. After this result is established, the convergence of the solution $y_h^{e_i}(\bomega,t,x)$ to $y^{e_i}(t,x)$ is a simple corollary that will be proven afterwards. 

\begin{lemma} \label{lem:conv}
Let $T > 0$ and suppose that $\mathbf{y}_0$, $\mathbf{y}_1$, and $\mathbf{u}$ satisfy the assumptions in Theorem \ref{thm:dyn2}. Then there exist a positive constant $C_T$ such that the solutions $w_{\pm}^{e_i}(t,x)$ to \eqref{eq:D} and the solutions $w_{h\pm}^{e_i}(\bomega,t,x)$ to \eqref{eq:RD} satisfy for all $e_i \in E$ and all $0 \leq t \leq T$
\begin{align*}
    \mathbb{E}[|w_{h+}^{e_i}(t)-w_{+}^{e_i}(t)|^2_{L^\infty(0,\ell_{e_i})}]&\leq C_T h(\mathrm{Lip}(\mathbf{y}_{0,x}) + \mathrm{Lip}(\mathbf{y}_{1}) + \mathrm{Lip}(\mathbf{u}))^2, \\
         \mathbb{E}[|w_{h-}^{e_i}(t)-w_{-}^{e_i}(t)|^2_{L^\infty(0,\ell_{e_i})}] &\leq C_T h(\mathrm{Lip}(\mathbf{y}_{0,x}) + \mathrm{Lip}(\mathbf{y}_{1}) + \mathrm{Lip}(\mathbf{u}))^2. 
\end{align*}
The constant $C_T$ depends on $T$ but is independent of the considered time grid $t_0, t_1, \ldots, t_K$, initial conditions $\mathbf{y}_0$ and $\mathbf{y}_1$, and the control $\mathbf{u}$. 
\end{lemma}

\begin{proof}
Write
$$
\delta = \min\left\{ \frac{\ell_{e_i}\pi_{e_i}}{c_{e_i}} \right\}_{e_i \in E} > 0,  
\qquad
\begin{bmatrix}
        e^{e_i}_{h+}(\bomega, t,x) \\
        e^{e_i}_{h-}(\bomega, t,x)
    \end{bmatrix} = \begin{bmatrix}
        w_{h+}^{e_i}(\bomega,t,x) - w_{+}^{e_i}(t,x) \\
        w_{h-}^{e_i}(\bomega,t,x) - w_{-}^{e_i}(t,x)
    \end{bmatrix},
    $$
Note that $\delta$ does not depend on the time grid used for the randomization. In particular, $\delta$ does not change when $h$ goes to zero. Also note that characteristics (both of the original and of the randomized system) require at least a time $\delta$ to travel from one endpoint of an edge to the other. This means that a characteristic can have at most one reflection in a time interval of length $\delta$. Furthermore, because the assumptions of Theorem \ref{thm:Lipschitz} are fulfilled, there exists a constant $\mathrm{Lip}(\mathbf{w})$ such that for all $e_i \in E$, all $x'', x' \in (0,\ell_{e_i})$, and all $t'', t' \in (0,T)$
\begin{subequations}
\label{eq:lem_Lipw}
\begin{align}
&|w^{e_i}_+(t'',x'') - w^{e_i}_+(t',x')| \leq \mathrm{Lip}(\mathbf{w}) (|x'' - x'| + \bar{c}|t'' - t'|),  \\
&|w^{e_i}_-(t'',x'') - w^{e_i}_-(t',x')| \leq \mathrm{Lip}(\mathbf{w}) (|x'' - x'| + \bar{c}|t'' - t'|), 
\end{align}
\end{subequations}
where $\bar{c} = \max_{e_i \in E} c_{e_i}$. 
Note that Theorem \ref{thm:Lipschitz} in particular shows that
\begin{equation}
\mathrm{Lip}(\mathbf{w}) \leq C_{\mathrm{Lip}}\left( \mathrm{Lip}(\mathbf{y}_{0,x}) + \mathrm{Lip}(\mathbf{y}_{1}) + \mathrm{Lip}(\mathbf{u}) \right) 
\label{eq:thm1_Lipw}
\end{equation}

The estimate is now proven by induction over the time instances $k\delta$. Note that $\delta$, and thus also the number of induction steps required to reach a certain time instant $t$, is \emph{not} dependent on the time grid used for randomization and does not increase when $h$ approaches zero. 

Our induction hypothesis is thus that there exists a constant $C_k$ independent of $h$ such that for all $e_i \in E$
\begin{equation}
\mathbb{E}[ |e^{e_i}_{h+}(k\delta) |_{L^\infty(0,\ell_{e_i})}^2 ] \leq C_k h (\mathrm{Lip}(\mathbf{w}))^2, \quad \mathbb{E}[ |e^{e_i}_{h-}(k\delta) |_{L^\infty(0,\ell_{e_i})}^2 ] \leq C_k h (\mathrm{Lip}(\mathbf{w}))^2,  \label{eq:inductionP}
\end{equation}
and we need to show that there exists a constant $C_{k+1}$ independent of $h$ such that for any $t\in [k\delta, (k+1)\delta]$ and $e_i \in E$
\begin{equation}\label{eq:toproofP}
     \mathbb{E}[ |e^{e_i}_{h+}(t) |_{L^\infty(0,\ell_{e_i})}^2 ] \leq C_{k+1} h (\mathrm{Lip}(\mathbf{w}))^2, \quad \mathbb{E}[ |e^{e_i}_{h-}(t) |_{L^\infty(0,\ell_{e_i})}^2 ] \leq C_{k+1} h (\mathrm{Lip}(\mathbf{w}))^2.
    \end{equation} 
Note that \eqref{eq:inductionP} holds for $k = 0$ because $w^{e_i}_{h\pm}(\bomega, 0, x) = w^{e_i}_\pm (0,x)$ for all $\bomega \in \mathcal{P}$.

The solutions $w_{+}^{e_i}(t,x)$ and $w_{h+}^{e_i}(\bomega,t,x)$ are traced back along characteristics. Each of the characteristics $s \mapsto (s,\xi_{+}^{e_i}(s; t,x))$ and $s \mapsto (s,\xi_{h+}^{e_i}(\bomega,s; t,x))$ can trace back to the line $\{ k\delta \} \times (0, \ell_{e_i})$ on which the induction hypothesis  \eqref{eq:inductionP} holds or to boundary $[k\delta,t] \times \{ 0 \}$. Therefore, four cases can appear when computing $e_{h+}^{e_i}(\bomega,t,x)$. 

In the first case, where both characteristcs trace back to $\{ k\delta \} \times (0, \ell_{e_i})$
\begin{align}
    e^{e_i}_{h+}(\bomega,t,x) &= w_{h+}^{e_i}(\bomega,k\delta, \xi^{e_i}_{h+}(\bomega,k\delta;t,x)) - w_+^{e_i}(k\delta, \xi^{e_i}_{+}(k\delta; t,x)) \nonumber \\
    &= e_{h+}^{e_i}(\bomega,k\delta, \xi^{e_i}_{h+}(\bomega,k\delta;t,x)) \nonumber \\
    &\qquad\qquad + w_+^{e_i}(k\delta, \xi^{e_i}_{h+}(\bomega,k\delta;t,x)) -w_+^{e_i}(k\delta, \xi^{e_i}_{+}(k\delta; t,x)). \label{eq:thm1_Step1a}
\end{align} 
Taking norms in \eqref{eq:thm1_Step1a}, it follows that
\begin{align}
    |e^{e_i}_{h+}(\bomega,t,x)| &\leq |e_{h+}^{e_i}(\bomega,k\delta)|_{L^\infty(0,\ell_{e_i})} + \mathrm{Lip}(\mathbf{w})|\xi^{e_i}_{h+}(\bomega,k\delta;t,x)) - \xi^{e_i}_{+}(k\delta; t,x))|. \label{eq:thm1_Step1b}
\end{align}
Note that the expected value of these two terms can be estimated by the induction hypothesis \eqref{eq:inductionP} and Lemma \ref{lem:conv_xi}. 

In the second case, where the randomized characteristic traces back to $\{ k\delta \} \times (0, \ell_{e_i})$ and the characteristic of the original system traces back to $[k\delta,t] \times \{ 0 \}$, 
\begin{align}
    e^{e_i}_{h+}(\bomega,t,x) &= w_{h+}^{e_i}(\bomega,k\delta, \xi^{e_i}_{h+}(\bomega,k\delta;t,x)) - w_+^{e_i}(t_{h+,\mathrm{in}}^{e_i}(t,x), 0) \nonumber \\
    &= e_{h+}^{e_i}(\bomega,k\delta, \xi^{e_i}_{h+}(\bomega,k\delta;t,x))+w_+^{e_i}(k\delta, \xi^{e_i}_{h+}(\bomega,k\delta;t,x)) - w_+^{e_i}(k\delta, 0) \nonumber \\
    &\qquad\qquad +w_+^{e_i}(k\delta, 0)-
    w_+^{e_i}(t_{h+,\mathrm{in}}^{e_i}(t,x), 0). 
\end{align}
Just as before, taking norms shows that
\begin{align}
    |e^{e_i}_{h+}&(\bomega,t,x)| \leq |e_{h+}^{e_i}(\bomega,k\delta)|_{L^\infty(0,\ell_{e_i})} + \mathrm{Lip}(\mathbf{w}) |\xi^{e_i}_{h+}(\bomega,k\delta;t,x)) - \xi^{e_i}_{+}(k\delta; t,x))| \nonumber \\
    & + \mathrm{Lip}(\mathbf{w}) \bar{c}|\max\{ t_{h+,\mathrm{in}}^{e_i}(\bomega,t,x), k\delta \} - \max\{ t_{+,\mathrm{in}}^{e_i}(t,x), k\delta \}|, \label{eq:thm1_Step2b}
\end{align}
where it has been used that $\xi^{e_i}_{+}(k\delta; t,x) < 0$ when the characteristic for the original system traces back to the boundary and that $t^{e_i}_{h+,\mathrm{in}}(\bomega,t,x) < k\delta$ when the randomized characteristic traces back to the line $\{ k\delta \} \times (0, \ell_{e_i})$. Note that the expected value of the three remaining terms can be estimated using the induction hypothesis \eqref{eq:inductionP}, Lemma \ref{lem:conv_xi}, and Lemma \ref{lem:conv_tini}. 

In the third case, where the randomized characteristic traces back to the boundary $[k\delta, t] \times \{ 0 \}$ and the characteristic of the original system to the line $\{k\delta \} \times (0, \ell_{e_i})$,
\begin{align}
    e^{e_i}_{h+}(\bomega,t,x) &= w_{h+}^{e_i}(\bomega, t_{h+,\mathrm{in}}^{e_i}(\bomega,t,x), 0) - w_+^{e_i}(k\delta, \xi^{e_i}_{+}(k\delta; t,x)) \nonumber \\
    &= e_{h+}^{e_i}(\bomega,t_{h+,\mathrm{in}}^{e_i}(\bomega,t,x), 0) + w_+^{e_i}(t_{h+,\mathrm{in}}^{e_i}(\bomega,t,x), 0) -w_+^{e_i}(k\delta, 0) \nonumber \\
    &\qquad\qquad + w_+^{e_i}(k\delta, 0) -w_+^{e_i}(k\delta, \xi^{e_i}_{+}(k\delta; t,x)).
\end{align}
Taking norms, it follows that 
\begin{align}
    |e^{e_i}_{h+}&(\bomega,t,x)| \leq |e_{h+}^{e_i}(\bomega,t^{e_i}_{h+,\mathrm{in}}(\bomega,t,x),0)| \nonumber \\
    & + \mathrm{Lip}(\mathbf{w}) \bar{c}|\max\{ t_{h+,\mathrm{in}}^{e_i}(\bomega,t,x), k\delta \} - \max\{ t_{+,\mathrm{in}}^{e_i}(t,x), k\delta \}| \nonumber \\
    & + \mathrm{Lip}(\mathbf{w}) |\xi^{e_i}_{h+}(\bomega,k\delta;t,x)) - \xi^{e_i}_{+}(k\delta; t,x))|,  \label{eq:thm1_Step3b}
\end{align}
where it has been used that $t^{e_i}_{+,\mathrm{in}}(t,x) \leq k\delta$ and $\xi^{e_i}_{h+}(\bomega,k\delta;t,x)) \leq 0$ in this case. Note that $|e_{h+}^{e_i}(\bomega,t^{e_i}_{h+,\mathrm{in}}(\bomega,t,x),0)|$ cannot be estimated by $|e^{e_i}_{h+}(k\delta)|_{L^\infty(0, \ell_{e_i})}$ because $t^{e_i}_{h+,\mathrm{in}}(\bomega,t,x) \geq k\delta$ in this case. This problem will also appear in the fourth case. After considering the fourth case, we will derive an additional estimate to deal with this term. 

In the fourth and final case the characteristics both trace back to the boundary $[k\delta, t] \times \{ 0 \}$ and
\begin{align}
    &e^{e_i}_{h+}(\bomega,t,x) = w_{h+}^{e_i}(\bomega, t_{h+,\mathrm{in}}^{e_i}(\bomega,t,x), 0) - w_+^{e_i}(t_{+,\mathrm{in}}^{e_i}(t,x), 0) \nonumber \\
    &\qquad = e_{h+}^{e_i}(\bomega,t_{h+,\mathrm{in}}^{e_i}(\bomega,t,x), 0) + w_+^{e_i}(t_{h+,\mathrm{in}}^{e_i}(\bomega,t,x), 0) -w_+^{e_i}(t_{+,\mathrm{in}}^{e_i}(t,x), 0),
\end{align}
so that after taking norms, it follows that
\begin{align}
    |e^{e_i}_{h+}&(\bomega,t,x)| \leq |e_{h+}^{e_i}(\bomega,t^{e_i}_{h+,\mathrm{in}}(\bomega,t,x),0)| \nonumber \\
    & + \mathrm{Lip}(\mathbf{w})\bar{c}|\max\{ t_{h+,\mathrm{in}}^{e_i}(\bomega,t,x), k\delta \} - \max\{ t_{+,\mathrm{in}}^{e_i}(t,x), k\delta \}|.  \label{eq:thm1_Step4b}
\end{align}

We now deal with the term $|e_{h+}^{e_i}(\bomega,t^{e_i}_{h+,\mathrm{in}}(\bomega,t,x),0)|$ that appears in the third and fourth case. 
In the following, write for brevity
\begin{equation}
 t^*_h(\bomega) = \max\{ t^{e_i}_{h+,\mathrm{in}}(\bomega,t,x), k\delta \}, \qquad 
 t^* = \max\{ t^{e_i}_{+,\mathrm{in}}(t,x), k\delta \}. 
\end{equation}
To find an estimate for $|e_{h+}^{e_i}(\bomega,t^*_h(\bomega),0)|$, let $v_j \in V$ be the node corresponding to the point $x=0$ on the edge $e_i \in E$. Subtracting the boundary conditions in \eqref{eq:D} from the boundary conditions in \eqref{eq:RD} shows that
\begin{equation}
   e_{h+}^{e_i}(\bomega,t^*_h(\bomega),0) = - e_{h,\mathrm{out}}^{e_i}(\bomega,t^*_h(\bomega),v_j) + \frac{2}{c^{v_j}_{\mathrm{tot}}}
    \sum_{e_\ell \in E(v_j)} c_{e_\ell}e^{e_\ell}_{h,\mathrm{out}}(\bomega, t^*_h(\bomega),v_j). \label{eq:thm1_Step5a}
\end{equation}
Because the choice of $\delta$ assures that $w^{e_\ell}_{h,\mathrm{out}}(\bomega, t^*_h(\bomega), v_j)$ and $w^{e_\ell}_{\mathrm{out}}(t^*_h(\bomega), v_j)$ trace back to the line $\{ k\delta \} \times (0, \ell_{e_\ell})$, estimates like the ones in the first case also hold for each $e^{e_\ell}_{h,\mathrm{out}}(\bomega, t^*_h(\bomega), v_j)$. In particular, $e^{e_\ell}_{h,\mathrm{out}}(\bomega,t^*_h(\bomega),v_j)$ can either be equal to outflow of a positive-going wave, i.e.\ $e^{e_\ell}_{h,\mathrm{out}}(\bomega,t^*_h(\bomega),v_j) = e^{e_\ell}_{h+}(\bomega,t^*_h(\bomega),\ell_{e_\ell})$, or the outflow of a negative-going wave, i.e.\ $e^{e_\ell}_{h,\mathrm{out}}(\bomega,t^*_h(\bomega),v_j) = e^{e_\ell}_{h-}(\bomega,t^*_h(\bomega),0)$. For a positive going wave, an estimate like \eqref{eq:thm1_Step1b} holds with $t = t^*_h(\bomega)$ and $x = \ell_{e_\ell}$ and for a negative going wave it can be shown similarly that
\begin{align}
	|e^{e_\ell}_{h,\mathrm{out}}(\bomega,t^*_h(\bomega),v_j)| &\leq |e_{h-}^{e_\ell}(\bomega,k\delta)|_{L^\infty(0,\ell_{e_i})} \nonumber \\
	&\qquad + \mathrm{Lip}(\mathbf{w})|\xi^{e_\ell}_{h-}(\bomega,k\delta;t^*_h(\bomega),0)) - \xi^{e_\ell}_{-}(k\delta; t^*_h(\bomega),0))|.
\end{align}
Inserting this estimate back into \eqref{eq:thm1_Step5a} now gives that
\begin{align}
&|e_{h+}^{e_i}(\bomega,t^*_h(\bomega),0)|\leq | e^{e_i}_{h-}(\bomega,k\delta)|_{L^\infty(0,\ell_{e_i})} \nonumber \\
&\qquad+ \frac{2}{c^{v_j}_{\mathrm{tot}}} \left( \sum_{e_\ell \in E^+(v_j)} c_{e_\ell}| e^{e_\ell}_{h+}(\bomega,k\delta)|_{L^\infty(0,\ell_{e_\ell})} + \sum_{e_\ell \in E^-(v_j)} c_{e_\ell}| e^{e_\ell}_{h-}(\bomega,k\delta)|_{L^\infty(0,\ell_{e_\ell})} \right) \nonumber  \\
&\qquad+ \mathrm{Lip}(\mathbf{w})|\xi^{e_i}_{h-}(\bomega,k\delta;t^*_h(\bomega),0)) - \xi^{e_i}_{-}(k\delta;t^*_h(\bomega),0))| \nonumber \\
	&\qquad+\frac{2\mathrm{Lip}(\mathbf{w})}{c^{v_j}_{\mathrm{tot}}} \sum_{e_\ell \in E^+(v_j)} c_{e_\ell}|\xi^{e_\ell}_{h+}(\bomega,k\delta;t^*_h(\bomega),\ell_{e_\ell})) - \xi^{e_\ell}_{+}(k\delta;t^*_h(\bomega),\ell_{e_\ell}))| \nonumber \\
	&\qquad+\frac{2\mathrm{Lip}(\mathbf{w})}{c^{v_j}_{\mathrm{tot}}} \sum_{e_\ell \in E^-(v_j)} c_{e_\ell}|\xi^{e_\ell}_{h-}(\bomega,k\delta;t^*_h(\bomega),0)) - \xi^{e_\ell}_{-}(k\delta;t^*_h(\bomega),0))|. \label{eq:thm1_Step5b} \\
	&\leq | e^{e_i}_{h-}(\bomega,k\delta)|_{L^\infty(0,\ell_{e_i})} \nonumber \\
&\qquad+ \frac{2}{c^{v_j}_{\mathrm{tot}}} \left( \sum_{e_\ell \in E^+(v_j)} c_{e_\ell}| e^{e_\ell}_{h+}(\bomega,k\delta)|_{L^\infty(0,\ell_{e_\ell})} + \sum_{e_\ell \in E^-(v_j)} c_{e_\ell}| e^{e_\ell}_{h-}(\bomega,k\delta)|_{L^\infty(0,\ell_{e_\ell})} \right) \nonumber  \\ 
	&\qquad+ \mathrm{Lip}(\mathbf{w})|\xi^{e_i}_{h-}(\bomega,k\delta;t^*,0)) - \xi^{e_i}_{-}(k\delta; t^*,0))| \nonumber \\ 
	&\qquad+\frac{2\mathrm{Lip}(\mathbf{w})}{c^{v_j}_{\mathrm{tot}}} \sum_{e_\ell \in E^+(v_j)} c_{e_\ell}|\xi^{e_\ell}_{h+}(\bomega,k\delta;t^*,0)) - \xi^{e_\ell}_{+}(k\delta; t^*,0))| \nonumber \\
	&\qquad+\frac{2\mathrm{Lip}(\mathbf{w})}{c^{v_j}_{\mathrm{tot}}} \sum_{e_\ell \in E^-(v_j)} c_{e_\ell}|\xi^{e_\ell}_{h-}(\bomega,k\delta;t^*,0)) - \xi^{e_\ell}_{-}(k\delta; t^*,0))| \nonumber \\
	&\qquad+3\mathrm{Lip}(\mathbf{w})\left(\bar{c}_h + \bar{c}\right)|t^*_h(\bomega) - t^*| \label{eq:thm1_Step5c}
\end{align}
where $E^+(v_j)$ denotes the set of all $e_\ell \in E(v_j)$ such that  $e^{e_\ell}_{h,\mathrm{out}}(\bomega,t,v_j) = e^{e_\ell}_{h+}(\bomega,t,\ell_{e_\ell})$, $E^-(v_j)$ denotes the set of all $e_\ell \in E(v_j)$ such that $e^{e_\ell}_{h,\mathrm{out}}(\bomega,t,v_j) = e^{e_\ell}_{h-}(\bomega,t,0)$. 
In the second inequality, it has been used that the characteristics $\xi_{h\pm}^{e_i}(\bomega,s;t,x)$ are Lipschitz in $t$ with Lipschitz constant $\bar{c}_h := \max_{e_i \in E}\{ c_{e_i}/\pi_{e_i} \}$ and the characteristics $\xi_{\pm}^{e_i}(s;t,x)$ are Lipschitz in $t$ with Lipschitz constant $\bar{c} := \max_{e_i \in E}\{ c_{e_i} \}$. Note that the shift to a deterministic time instant $t^*$ in the second inequality is needed to apply the error estimate from Lemma \ref{lem:conv_xi}. 

Combining \eqref{eq:thm1_Step1b}, \eqref{eq:thm1_Step2b}, \eqref{eq:thm1_Step3b}, \eqref{eq:thm1_Step4b}, and \eqref{eq:thm1_Step5b} shows that in any of the four cases
\begin{align}
&|e^{e_i}_{h+}(\bomega,t,x)| \leq | e^{e_i}_{h-}(\bomega,k\delta)|_{L^\infty(0,\ell_{e_i})} \nonumber \\
&\quad+ \frac{2}{c^{v_j}_{\mathrm{tot}}} \left( \sum_{e_\ell \in E^+(v_j)} c_{e_\ell}| e^{e_\ell}_{h+}(\bomega,k\delta)|_{L^\infty(0,\ell_{e_\ell})} + \sum_{e_\ell \in E^-(v_j)} c_{e_\ell}| e^{e_\ell}_{h-}(\bomega,k\delta)|_{L^\infty(0,\ell_{e_\ell})} \right) \nonumber  \\ 
	&\qquad+ \mathrm{Lip}(\mathbf{w})|\xi^{e_i}_{h+}(\bomega,k\delta;t)) - \xi^{e_i}_{+}(k\delta; t))|_{L^\infty(0,\ell_{e_i})} \nonumber \\
	&\qquad+ \mathrm{Lip}(\mathbf{w})|\xi^{e_i}_{h-}(\bomega,k\delta;t^*)) - \xi^{e_i}_{-}(k\delta; t^*))|_{L^\infty(0,\ell_{e_i})} \nonumber \\ 
	&\qquad+\frac{2\mathrm{Lip}(\mathbf{w})}{c^{v_j}_{\mathrm{tot}}} \sum_{e_\ell \in E^+(v_j)} c_{e_\ell}|\xi^{e_\ell}_{h+}(\bomega,k\delta;t^*)) - \xi^{e_\ell}_{+}(k\delta; t^*))|_{L^\infty(0,\ell_{e_i})} \nonumber \\
	&\qquad+\frac{2\mathrm{Lip}(\mathbf{w})}{c^{v_j}_{\mathrm{tot}}} \sum_{e_\ell \in E^-(v_j)} c_{e_\ell}|\xi^{e_\ell}_{h-}(\bomega,k\delta;t^*)) - \xi^{e_\ell}_{-}(k\delta; t^*))|_{L^\infty(0,\ell_{e_i})} \nonumber \\
	&\qquad+\mathrm{Lip}(\mathbf{w})\left(3\bar{c}_h + 4\bar{c}\right)|\max\{ t_{h+,\mathrm{in}}^{e_i}(\bomega,t), k\delta \} - \max\{ t_{+,\mathrm{in}}^{e_i}(t), k\delta \}|_{L^\infty(0,\ell_{e_i})}. \label{eq:thm1_Step6}
\end{align}
Note that the LHS can be replaced by $|e_{h+}(\bomega,t)|_{L^\infty(0,\ell_{e_i})}$ because the RHS of \eqref{eq:thm1_Step6} does not depend on $x$. Taking the expected value using $\sqrt{\mathbb{E}[(X+Y)^2]} \leq \sqrt{\mathbb{E}[X^2]} + \sqrt{\mathbb{E}[Y^2]}$, the induction hypothesis \eqref{eq:inductionP}, \eqref{eq:lem_xi1} in Lemma \ref{lem:conv_xi}, and Lemma \ref{lem:conv_tini}, it follows that

\begin{multline}
\sqrt{\mathbb{E}[|e^{e_i}_{h+}(t)|^2_{L^\infty(0,\ell_{e_i})}]} 
\leq 3 \mathrm{Lip}(\mathbf{w})\sqrt{C_k h} + 4 \mathrm{Lip}(\mathbf{w})\sqrt{h(t-k\delta)\mathrm{Var}[c_h]} \\
    + \mathrm{Lip}(\mathbf{w})(3\bar{c}_h + 4\bar{c}) \sqrt{C_2 h(t-k\delta)}.  \label{eq:thm1_Step7}
\end{multline}
It is now easy to verify that the first inequality in \eqref{eq:toproofP} holds for a certain constant $C_{k+1}$. 
The second inequality in \eqref{eq:toproofP} for $\mathbb{E}[|e^{e_i}_{h+}(t)|^2_{L^\infty(0,\ell_{e_i})}]$ can be proven completely analogous. The statement of the theorem then follows after using \eqref{eq:thm1_Lipw}. 
\end{proof}

The estimate for the Riemann variables in Lemma \ref{lem:conv} implies a similar estimate on the solutions to the wave equation
\begin{align}
y_h^{e_i}(\bomega,t,x) &= y_0^{e_i}(x) + \int_0^t \frac{ w_{h+}^{e_i}(\bomega,s,x) + w_{h-}^{e_i}(\bomega,s,x) }{2} \ \mathrm{d}s, \label{eq:yhprime} \\
y^{e_i}(\bomega,t,x) &= y_0^{e_i}(x) + \int_0^t \frac{ w_{+}^{e_i}(\bomega,s,x) + w_{-}^{e_i}(\bomega,s,x) }{2} \ \mathrm{d}s.
\label{eq:yprime}
\end{align}

\begin{corollary} \label{corr:yh_y}
Under the assumptions of Theorem \ref{thm:dyn2}, there also exists a constant $C_T$ independent of the considered time grid such that for all $e_i \in E$ and all $0 \leq t \leq T$
\begin{equation}
\mathbb{E}[|y_{h}^{e_i}(t)-y^{e_i}(t)|^2_{L^\infty(0,\ell_{e_i})}]\leq C_T h(\mathrm{Lip}(\mathbf{y}_{0,x}) + \mathrm{Lip}(\mathbf{y}_{1}) + \max_{\bomega \in \mathcal{P}} \mathrm{Lip}(\mathbf{u}_h(\bomega)))^2.
\end{equation}
\end{corollary}

\begin{proof}
By subtracting \eqref{eq:yprime} from \eqref{eq:yhprime}, it follows that
\begin{align}
y_h^{e_i}(\bomega,t,x) - y^{e_i}(\bomega,t,x) = \int_0^t \frac{ e^{e_i}_{h+}(\bomega,s,x) + e^{e_i}_{h-}(\bomega,s,x)}{2} \ \mathrm{d}s, 
\end{align}
where $e^{e_i}_{h+}(\bomega,s,x) = w^{e_i}_{h+}(\bomega,s,x) - e^{e_i}_{+}(\bomega,s,x)$ and $e^{e_i}_{h-}(\bomega,s,x) = w^{e_i}_{h-}(\bomega,s,x) - e^{e_i}_{-}(\bomega,s,x)$ as in the proof of Lemma \ref{lem:conv}. 
Taking the supremum over $x \in (0,\ell_{e_i})$ and applying Cauchy-Schwarz in time, it follows that 
\begin{align}
|y_h^{e_i}(\bomega,t) - y^{e_i}(\bomega,t)|_{L^\infty(0,\ell_{e_i})}^2 \leq \frac{t}{2} \int_0^t \left( |e^{e_i}_{h+}(\bomega,s)|_{L^\infty(0, \ell_{e_i})}^2 + |e^{e_i}_{h-}(\bomega,s)|_{L^\infty(0,\ell_{e_i})}^2 \right) \ \mathrm{d}s. 
\end{align}
The result now follows directly by taking the expectation and using the estimates from Lemma \ref{lem:conv}. 
\end{proof}

Lemma \ref{lem:conv} and Corollary \ref{corr:yh_y} together prove Theorem \ref{thm:dyn2}. 

\subsection{A Convergence Result for Random Controls}
For our analysis of the randomized optimal control problem in the next section, a second convergence result similar to Theorem \ref{thm:dyn2} is needed in which the control $\mathbf{u}$ is replaced by a random control $\mathbf{u}^*_h(\bomega,t)$. In particular, consider the solutions to 
\begin{equation} \left\{
\begin{aligned}
& w_{-,t}^{e_i}(\bomega,t,x) - c_{e_i} w_{-,x}^{e_i}(\bomega,t,x) = 0, \\
& w_{+,t}^{e_i}(\bomega,t,x) + c_{e_i} w_{+,x}^{e_i}(\bomega,t,x) = 0, \\
    & w^{e_i}_{\mathrm{in}}(\bomega,t,v_j) = - w^{e_i}_{\mathrm{out}}(\bomega,t,v_j) + \frac{2}{c_{\mathrm{tot}}^{v_j}}\left( \sum_{e_\ell \in E(v_j)}  c_{e_\ell}w_{\mathrm{out}}^{e_\ell}(\bomega,t,v_j) -\bar{u}_h^{v_j}(\bomega,t)\right), \\
    & w^{e_i}_-(\bomega,0,x) =  y^{e_i}_1 (x) + c_{e_i} y_{0,x}^{e_i}(x), 
    \qquad w^{e_i}_+(\bomega,0,x) =  y^{e_i}_1 (x) -c_{e_i} y_{0,x}^{e_i}(x), 
    \end{aligned}\right. 
    \tag{D'} \label{eq:Dprime}
\end{equation}
and 
\begin{equation} \left\{
\begin{aligned}
& w_{h-,t}^{e_i}(\bomega,t,x) - c_{h,e_i}(\bomega,t) w_{h-,x}^{e_i}(\bomega,t,x) = 0, \\
& w_{h+,t}^{e_i}(\bomega,t,x) + c_{h,e_i}(\bomega,t) w_{h+,x}^{e_i}(\bomega,t,x) = 0, \\
    & w^{e_i}_{h,\mathrm{in}}(\bomega,t,v_j) = - w^{e_i}_{h,\mathrm{out}}(\bomega,t,v_j) + \frac{2}{c_{\mathrm{tot}}^{v_j}} \left(\sum_{e_\ell \in E(v_j)} c_{e_\ell} w^{e_\ell}_{h,\mathrm{out}}(\bomega,t,v_j) -  \bar{u}_h^{v_j}(\bomega,t) \right), \\
    & w_{h-}^{e_i}(\bomega,0,x) = y^{e_i}_1 (x) + c_{e_i} y_{0,x}^{e_i}(x), \qquad  w_{h+}^{e_i}(\bomega,0,x) = y^{e_i}_1 (x) - c_{e_i} y_{0,x}^{e_i}(x). \\
    \end{aligned}\right. 
    \tag{RD'} \label{eq:RDprime}
\end{equation}
Note that $w_{\pm}^{e_i}(\bomega,t,x)$ now also depends on $\bomega$ because the control $\mathbf{u}_h(\bomega,t)$ depends on $\bomega$, but it still follows deterministic dynamics. The following error estimate will be used to prove the convergence of the optimal control in the next section. 

\begin{lemma} \label{lem:conv_uh}
Let $T > 0$ and suppose that $\mathbf{y}_0$, $\mathbf{y}_1$, and $\mathbf{u}_h(\bomega)$ are such that 1) they satisfy the boundary conditions in \eqref{eq:Dprime}, 2) $y_{0,x}^{e_i}(x)$ and $y_1^{e_i}(x)$ are Lipschitz on $(0, \ell_{e_i})$ for each $e_i \in E$, and 3) $\mathbf{u}_h(\bomega, t)$ is Lipschitz in $t$ on $(0,T)$ for all $\bomega \in \mathcal{P}$, 
then there exist a positive constant $C_T$ such that the solutions $w_{\pm}^{e_i}(\bomega,t,x)$ to \eqref{eq:Dprime} and the solutions $w_{h\pm}^{e_i}(\bomega,t,x)$ to \eqref{eq:RDprime} satisfy for all $e_i \in E$ and all $0 \leq t \leq T$
\begin{align*}
    \mathbb{E}[|w_{h+}^{e_i}(t)-w_{+}^{e_i}(t)|_{L^\infty(0,\ell_{e_i})}]&\leq C_T \sqrt{h}\sqrt{\mathbb{E}[(\mathrm{Lip}(\mathbf{y}_{0,x}) + \mathrm{Lip}(\mathbf{y}_{1}) + \mathrm{Lip}(\mathbf{u}_h))^2]}, \\
         \mathbb{E}[|w_{h-}^{e_i}(t)-w_{-}^{e_i}(t)|_{L^\infty(0,\ell_{e_i})}] &\leq C_T \sqrt{h}\sqrt{\mathbb{E}[(\mathrm{Lip}(\mathbf{y}_{0,x}) + \mathrm{Lip}(\mathbf{y}_{1}) + \mathrm{Lip}(\mathbf{u}_h))^2]}. 
\end{align*}
The constant $C_T$ depends on $T$ but is further independent of the considered time grid $t_0, t_1, \ldots, t_K$, initial conditions $\mathbf{y}_0$ and $\mathbf{y}_1$, and the control $\mathbf{u}_h(\bomega)$. 
\end{lemma}

\begin{proof}
The proof follows essentially the same steps as the proof of Lemma \ref{lem:conv}, with some small adaptations. Naturally, the induction hypothesis is replaced by
\begin{equation}
\begin{aligned}
\mathbb{E}[ |e^{e_i}_{h+}(k\delta) |_{L^\infty(0,\ell_{e_i})} ] \leq C_k \sqrt{h} \sqrt{\mathbb{E}[(\mathrm{Lip}(\mathbf{w}))^2]}, \\ \mathbb{E}[ |e^{e_i}_{h-}(k\delta) |_{L^\infty(0,\ell_{e_i})} ] \leq C_k \sqrt{h} \sqrt{\mathbb{E}[(\mathrm{Lip}(\mathbf{w}))^2]}, 
\end{aligned} \label{eq:inductionP2omega}
\end{equation}

Furthermore, applying Theorem \ref{thm:Lipschitz} for each $\bomega \in \mathcal{P}$, there exists a constant $\mathrm{Lip}(\mathbf{w}(\bomega))$ such that for all $e_i \in E$, $x'', x' \in (0, \ell_{e_i})$, and $t'', t' \in [0,T]$
\begin{subequations}
\begin{align}
&|w^{e_i}_+(\bomega,t'',x'') - w^{e_i}_+(\bomega,t',x')| \leq \mathrm{Lip}(\mathbf{w}(\bomega)) (|x'' - x'| + \bar{c}|t'' - t'|),  \\
&|w^{e_i}_-(\bomega,t'',x'') - w^{e_i}_-(\bomega,t',x')| \leq \mathrm{Lip}(\mathbf{w}(\bomega)) (|x'' - x'| + \bar{c}|t'' - t'|), 
\end{align}
\end{subequations}
which satisfies
\begin{equation}
\mathrm{Lip}(\mathbf{w}(\bomega)) \leq C_{\mathrm{Lip}}\left( \mathrm{Lip}(\mathbf{y}_{0,x}) + \mathrm{Lip}(\mathbf{y}_{1}) + \mathrm{Lip}(\mathbf{u}_h(\bomega)) \right). 
\label{eq:thm1_Lipw_omega}
\end{equation}

Following exactly the same procedure as in the proof of Lemma \ref{lem:conv}, an estimate of the form \eqref{eq:thm1_Step6} in which $\mathrm{Lip}(\mathbf{w})$ is replaced by $\mathrm{Lip}(\mathbf{w}(\bomega))$ can be obtained. In particular,
\begin{align}
&|e^{e_i}_{h+}(\bomega,t,x)| \leq | e^{e_i}_{h-}(\bomega,k\delta)|_{L^\infty(0,\ell_{e_i})} \nonumber \\
&\quad+ \frac{2}{c^{v_j}_{\mathrm{tot}}} \left( \sum_{e_\ell \in E^+(v_j)} c_{e_\ell}| e^{e_\ell}_{h+}(\bomega,k\delta)|_{L^\infty(0,\ell_{e_\ell})} + \sum_{e_\ell \in E^-(v_j)} c_{e_\ell}| e^{e_\ell}_{h-}(\bomega,k\delta)|_{L^\infty(0,\ell_{e_\ell})} \right) \nonumber  \\ 
	&\qquad+ \mathrm{Lip}(\mathbf{w}(\bomega))|\xi^{e_i}_{h+}(\bomega,k\delta;t)) - \xi^{e_i}_{+}(k\delta; t))|_{L^\infty(0,\ell_{e_i})} \nonumber \\
	&\qquad+ \mathrm{Lip}(\mathbf{w}(\bomega))|\xi^{e_i}_{h-}(\bomega,k\delta;t^*)) - \xi^{e_i}_{-}(k\delta; t^*))|_{L^\infty(0,\ell_{e_i})} \nonumber \\ 
	&\qquad+\frac{2\mathrm{Lip}(\mathbf{w}(\bomega))}{c^{v_j}_{\mathrm{tot}}} \sum_{e_\ell \in E^+(v_j)} c_{e_\ell}|\xi^{e_\ell}_{h+}(\bomega,k\delta;t^*)) - \xi^{e_\ell}_{+}(k\delta; t^*))|_{L^\infty(0,\ell_{e_i})} \nonumber \\
	&\qquad+\frac{2\mathrm{Lip}(\mathbf{w}(\bomega))}{c^{v_j}_{\mathrm{tot}}} \sum_{e_\ell \in E^-(v_j)} c_{e_\ell}|\xi^{e_\ell}_{h-}(\bomega,k\delta;t^*)) - \xi^{e_\ell}_{-}(k\delta; t^*))|_{L^\infty(0,\ell_{e_i})} \nonumber \\
	&\qquad+\mathrm{Lip}(\mathbf{w}(\bomega))\left(3\bar{c}_h + 4\bar{c}\right)|\max\{ t_{h+,\mathrm{in}}^{e_i}(\bomega,t), k\delta \} - \max\{ t_{+,\mathrm{in}}^{e_i}(t), k\delta \}|_{L^\infty(0,\ell_{e_i})}. \label{eq:thm1_Step6_omega}
\end{align}
Now taking directly the expectation using that $\mathbb{E}[XY] \leq \sqrt{\mathbb{E}[X^2]} \sqrt{\mathbb{E}[Y^2]}$ for the terms involving $\mathrm{Lip}(\mathbf{w}(\bomega))$, it follows that
\begin{multline}
\mathbb{E}[|e^{e_i}_{h+}(t)|_{L^\infty(0,\ell_{e_i})}]
\leq 3 C_k \sqrt{h} \sqrt{\mathbb{E}[(\mathrm{Lip}(\mathbf{w}))^2]} + 4 \sqrt{\mathbb{E}[(\mathrm{Lip}(\mathbf{w}))^2]}\sqrt{h(t-k\delta)\mathrm{Var}[c_h]} \\
    + \sqrt{\mathbb{E}[(\mathrm{Lip}(\mathbf{w}))^2]}(3\bar{c}_h + 4\bar{c}) \sqrt{C_2 h(t-k\delta)},  \label{eq:thm1_Step8}
\end{multline}
where the new induction hypothesis \eqref{eq:inductionP2omega}, and Lemmas \ref{lem:conv_xi} and \ref{lem:conv_tini} have been used. 
Because $\mathrm{Lip}(\mathbf{w}(\bomega))$ can be bounded as in \eqref{eq:thm1_Lipw_omega}, this estimate (and the analogue result for $e_{h-}^{e_i}(\bomega,t,x)$) now easily implies the estimates in Lemma \ref{lem:conv_uh}. 
\end{proof}

\section{Convergence in Optimal Controls \label{sec:control}}
This section contains the proof of Theorem \ref{thm:contr}. 

In this section, $C_T$ denotes a constant depending on the time horizon $T$ but independent of the chosen time grid and $h$ that may vary from line to line. 

The proof relies the strict convexity of the functional $J_h(\bomega,\mathbf{u})$ in \eqref{eq:OCPh}, from which the following classical observations follow, see, for example, \cite{minoux1986}. Consider the cost functional $J_h(\bomega, \mathbf{u})$ in \eqref{eq:OCPh} and note that this cost functional is strictly convex in the sense that for any controls $\mathbf{u}, \mathbf{v} \in U_{ad}$ and  $\theta \in [0,1]$
\begin{equation}
J_h(\bomega,(1-\theta)\mathbf{u} + \theta\mathbf{v}) \leq (1-\theta)J_h(\bomega, \mathbf{u}) + \theta J_h(\bomega, \mathbf{v}) - \alpha \theta (1-\theta) |\mathbf{u} - \mathbf{v}|^2_{H^2(0,T)}. 
\end{equation} 
Subtracting $J_h(\bomega,\mathbf{u})$ at both sides, dividing by $\theta$, and then taking the limit $\theta \rightarrow 0$, it follows that
\begin{equation}
J_h(\bomega,\mathbf{v}) \geq J_h(\bomega,\mathbf{u}) + \delta J_h(\bomega,\mathbf{u}; \mathbf{v} - \mathbf{u}) + \alpha  | \mathbf{u}-\mathbf{v} |^2_{H^2(0,T)}, \label{eq:convex_diff}
\end{equation}
where $\delta J_h(\bomega, \mathbf{u}; \mathbf{v})$ is the G\^ateaux-derivative of $J_h(\bomega, \cdot)$ at the point $\mathbf{u}$ in direction $\mathbf{v}$
\begin{equation}
\delta J_h(\bomega, \mathbf{u}; \mathbf{v}) = \lim_{\varepsilon \rightarrow 0} \frac{J_h(\bomega, \mathbf{u} + \varepsilon \mathbf{v}) - J_h(\bomega, \mathbf{u})}{\varepsilon}. \label{eq:gateaux}
\end{equation}

Recall that $\mathbf{u}^*(t)$ is the solution to the OCP in \eqref{eq:OCP} and that $\mathbf{u}^*_h(\bomega,t)$ the solution to the randomized OCP in \eqref{eq:OCPh}. Taking $\mathbf{u} = \mathbf{u}^*$ and $\mathbf{v} = \mathbf{u}^*_h(\bomega)$ in \eqref{eq:convex_diff} shows that
\begin{equation}
J_h(\bomega,\mathbf{u}_h^*(\bomega)) \geq J_h(\bomega,\mathbf{u}^*) + \delta J_h(\bomega,\mathbf{u}^*; \mathbf{u}^*_h(\bomega) - \mathbf{u}^*) + \alpha  | \mathbf{u}^*_h(\bomega)-\mathbf{u}^* |^2_{H^2(0,T)}, \label{eq:convex_diff2}
\end{equation}
Observe that $J_h(\bomega, \mathbf{u}^*_h(\bomega)) \leq J_h(\bomega,\mathbf{u}^*)$
because $\mathbf{u}^*_h(\bomega)$ is the minimizer of $J_h(\bomega, \cdot)$. Therefore, 
\begin{equation}
\delta J_h(\bomega,\mathbf{u}^*; \mathbf{u}^*_h(\bomega) - \mathbf{u}^*) + \alpha  | \mathbf{u}^*_h(\bomega)-\mathbf{u}^* |^2_{H^2(0,T)} \leq 
J_h(\bomega,\mathbf{u}_h^*(\bomega)) -  J_h(\bomega,\mathbf{u}^*) \leq 0, \nonumber  
\end{equation}
which in turn implies that
\begin{equation} 
\alpha  | \mathbf{u}^*_h(\bomega)-\mathbf{u}^* |^2_{H^2(0,T)} \leq 
|\delta J_h(\bomega,\mathbf{u}^*; \mathbf{u}^*_h(\bomega) - \mathbf{u}^*)|.  \label{eq:convex_diff4}
\end{equation}
In order to prove Theorem \ref{thm:contr}, it thus suffices to bound the expected value of $|\delta J_h(\bomega,\mathbf{u}^*; \mathbf{u}^*_h(\bomega) - \mathbf{u}^*)|$. Such a bound is provided by the following lemma. 

\begin{lemma} Consider a random perturbation $\mathbf{v}_h(\bomega,t) \in H^2(0,T; \mathbb{R}^{|V_C|})$ of the optimal control $\mathbf{u}^*(t)$ for \eqref{eq:OCP}. Then there exists a constant $C$, independent of $h$, such that \label{lem:contr}
\begin{equation}
\mathbb{E}[|\delta J_h(\mathbf{u}^*; \mathbf{v}_h)|] \leq C \sqrt{h} \sqrt{\mathbb{E}[|\mathbf{v}_h|_{H^2(0,T)}^2]}. \label{eq:lem_contr}
\end{equation}
\end{lemma}

\begin{proof}
Because $\mathbf{u}^*$ is the minimizer of $J(\cdot)$, $\delta J(\mathbf{u}^*; \mathbf{v}_h(\bomega)) = 0$ for any $\mathbf{v}_h(\bomega) \in \mathcal{U}^1$. Therefore,
\begin{equation}
|\delta J_h(\bomega,\mathbf{u}^*; \mathbf{v}^*_h(\bomega))| = |\delta J_h(\bomega,\mathbf{u}^*; \mathbf{v}^*_h(\bomega)) - \underbrace{\delta J(\mathbf{u}^*; \mathbf{v}^*_h(\bomega))}_{=0}| \label{eq:lem_contr_step1}
\end{equation} 
Now let $\mathbf{y}^*(t,x)$ and $\mathbf{y}^*_h(\bomega,t,x)$ denote the solutions to \eqref{eq:D} and \eqref{eq:RD} resulting from the control $\mathbf{u}^*(t)$ and the initial conditions $(\mathbf{y}_0, \mathbf{y}_1)$, respectively, and let $\mathbf{z}(\bomega,t,x)$ and $\mathbf{z}_h(\bomega,t,x)$ be the solutions to \eqref{eq:Dprime} and \eqref{eq:RDprime} resulting from a random perturbation of the control $\mathbf{v}^*_h(\bomega,t)$ and zero initial conditions, respectively.

Taking the G\^ateaux-derivative of $J_h(\bomega, \cdot)$ and $J(\cdot)$ according to \eqref{eq:gateaux} shows that
\begin{align}
    \delta J_h(\bomega, \mathbf u^*; \mathbf{v}^*_h(\bomega)) &= \langle \mathbf y_h^*(\bomega) - \mathbf y_d, \mathbf{z}_h(\bomega) \rangle_{L^2(Q)} +  \alpha \langle \mathbf{u}, \mathbf{v}_h(\bomega) \rangle_{H^2(0,T)}, \nonumber \\
	\delta J(\mathbf u^*; \mathbf{v}_h(\bomega)) &= \langle \mathbf y^* - \mathbf y_d, \mathbf{z}(\bomega)  \rangle_{L^2(Q)} + \alpha \langle \mathbf{u}, \mathbf{v}_h(\bomega) \rangle_{H^2(0,T)}.  \nonumber
\end{align}
Subtracting these two expressions shows that
\begin{align}
\delta J_h(\bomega,\mathbf{u}^*; \mathbf{v}_h(\bomega)) &- \delta J(\mathbf{u}^*; \mathbf{v}_h(\bomega) ) \nonumber \\
&= \langle \mathbf y_h^*(\bomega) - \mathbf y_d, \mathbf{z}_h(\bomega) \rangle_{L^2(Q)} - \langle \mathbf y^* - \mathbf y_d, \mathbf{z}(\bomega)  \rangle_{L^2(Q)} \nonumber \\
&= \langle \mathbf y_h^*(\bomega) - \mathbf y^*, \mathbf{z}_h(\bomega) \rangle_{L^2(Q)} - \langle \mathbf y^* - \mathbf y_d, \mathbf{z}_h(\bomega) - \mathbf{z}(\bomega)  \rangle_{L^2(Q)}.   \label{eq:lem_contr_step2}
\end{align}
Taking the absolute value and using Cauchy Schwarz in $L^2(Q)$, it follows that
\begin{align}
& | \delta J_h(\bomega,\mathbf{u}^*;\mathbf{v}^*_h(\bomega))| = | \delta J_h(\bomega,\mathbf{u}^*; \mathbf{v}^*_h(\bomega)) - \delta J(\mathbf{u}^*; \mathbf{v}^*_h(\bomega) )| \nonumber \\
&\qquad \leq |\mathbf y_h^*(\bomega) - \mathbf y^*|_{L^2(Q)} | \mathbf{z}_h(\bomega) |_{L^2(Q)} + | \mathbf y^* - \mathbf y_d|_{L^2(Q)} |\mathbf{z}_h(\bomega) - \mathbf{z}(\bomega)|_{L^2(Q)}. \label{eq:lem_contr_step3}
\end{align}
Now taking the expectation results in
\begin{multline}
\mathbb{E}[ | \delta J_h(\mathbf{u}^*;\mathbf{v}^*_h)|] \leq \sqrt{\mathbb{E}[|\mathbf y_h^* - \mathbf y^*|_{L^2(Q)} ^2]} \sqrt{\mathbb{E}[| \mathbf{z}_h |_{L^2(Q)}^2]} \\ + | \mathbf y^* - \mathbf y_d|_{L^2(Q)} \mathbb{E}[ |\mathbf{z}_h - \mathbf{z}|_{L^2(Q)}]. \label{eq:lem_contr_step4}
\end{multline}
By Lemma \ref{lem:conv}, there exists a constant $C_T$ (depending on the initial conditions and $\mathrm{Lip}(\mathbf{u}^*)$) such that
\begin{equation}
\mathbb{E}[|\mathbf y_h^* - \mathbf y^*|_{L^2(Q)} ^2] \leq C h. 
\end{equation}
By Theorem \ref{thm. wellposednessh}, there exists a constant $C_T$ such that $|\mathbf{z}_h(\bomega)|_{L^2(Q)} ^2 \leq C_T |\mathbf{v}_h(\bomega)|_{L^2(0,T)}^2$, which implies that also
\begin{equation}
\mathbb{E}[|\mathbf{z}_h|_{L^2(Q)}^2] \leq C_T \mathbb{E}[|\mathbf{v}_h|_{L^2(0,T)}^2] \leq C_T \mathbb{E}[|\mathbf{v}_h|_{H^2(0,T)}^2]. 
\end{equation}
By Theorem \ref{thm. wellposedness}, there also exists a constant $C_T$ such that $|\mathbf{y}^* - \mathbf{y}_d|_{L^2(Q)}^2 \leq C_T$.  \\
By Lemma \ref{lem:conv_uh}, there exists a constant $C_T$ such that
\begin{equation}
\mathbb{E}[|\mathbf{z}_h - \mathbf{z}|_{L^2(Q)}] \leq C_T \sqrt{h} \sqrt{\mathbb{E}[\mathrm{Lip}(\mathbf{v}_h)^2]} \leq C_T \sqrt{h} \sqrt{\mathbb{E}[|\mathbf{v}_h|_{H^2(0,T)}^2]}. 
\end{equation}

By inserting these four observations back into \eqref{eq:lem_contr_step4}, the result follows. 
\end{proof}

With this result, the proof of Theorem \ref{thm:contr} can now be completed as follows. 

\begin{proof}[Proof of Theorem \ref{thm:contr}]

Taking the expected value in \eqref{eq:convex_diff4} and inserting the result from Lemma \ref{lem:contr} with $\mathbf{v}_h(\bomega,t) = \mathbf{u}^*_h(\bomega,t) - \mathbf{u}^*(t)$, it follows that
\begin{equation}
\alpha \mathbb{E}[|\mathbf{u}^*_h - \mathbf{u}^*|^2_{H^2(0,T)}] \leq 
\mathbb{E}[|\delta J_h(\mathbf{u}^*; \mathbf{u}^*_h - \mathbf{u}^*)|] \leq C_T \sqrt{h} \sqrt{\mathbb{E}[|\mathbf{u}^*_h - \mathbf{u}^*|_{H^2(0,T)}^2]}. 
\end{equation}
Dividing by $\sqrt{\mathbb{E}[|\mathbf{u}^*_h - \mathbf{u}^*|_{H^2(0,T)}^2]}$ and squaring the result,  it follows that
\begin{equation}
\alpha^2 \mathbb{E}[|\mathbf{u}^*_h - \mathbf{u}^*|^2_{H^2(0,T)}] \leq C_T^2 h. 
\end{equation}
This is the result from Theorem \ref{thm:contr}. 
\end{proof}

\section{Numerical Examples \label{sec:examples}}

In this section, we validate the two main results of this paper (Theorem \ref{thm:dyn2} and Theorem \ref{thm:contr}) for two numerical examples. The first example considers the diamond network from Example \ref{example:diamond_graph} and the second example considers a larger network with more than 40 edges from the Gaslib library \cite{gaslib2017, gaslibreport2017}. 

\subsection{The Diamond Network \label{ssec:examples_diamond}}
Consider the diamond network from Example \ref{example:diamond_graph} in Figure \ref{fig:diamond-graph}. As demonstrated in Example \ref{example:diamond_graph_contd}, this graph can be decomposed into the three tree-shaped networks in Figure \ref{fig:ram_diamond}. As suggested in Figure \ref{fig:diamond-graph}, we take $\ell_{e_1} = \ell_{e_4} = \ell_{e_7} = \sqrt{2}$ and $\ell_{e_2} = \ell_{e_3} = \ell_{e_5} = \ell_{e_6} = 1$ and take $c_{e_i} = 1$ on all edges $e_i \in E$. The randomized velocity field $c_{h,e_i}(\bomega,t)$ can then be constructed as described in Example \ref{example:diamond_graph_contd}. We make a spatial discretiziation of \eqref{eq:D} and \eqref{eq:RD} by the first-order upwind scheme. The spatial grid is uniform on each edge. We use 30 spatial grid points on the edges of length $\sqrt{2}$ and 21 grid points on the edges of length 1, which assures that the  grid spacing $\Delta x \leq 0.05$ on all edges. Time is discretized with the backward Euler scheme with a uniform step size $h$. For the dynamics \eqref{eq:RD} one out of the four subgraphs in Figure \ref{fig:ram_diamond} is chosen for each discretization step (each with a probability $p_\omega = \tfrac{1}{4}$). The control is applied at node $v_1$, i.e.\ $V_C = \{ v_1 \}$. 

To validate the convergence of the forward dynamics from Theorem \ref{thm:dyn2}, the control is fixed as $u(t) = \sin(\pi t)$. The solutions $w_\pm(t,x)$ and $w_{h\pm}(\bomega,t,x)$ to \eqref{eq:D} and \eqref{eq:RD} are computed using the procedure described above starting from zero initial conditions over a time interval of length $T = 5$. The results are summarized in Table \ref{tab:Diamond_dyn}. To get an idea of the variance in the obtained results, the solutions to \eqref{eq:D} and \eqref{eq:RD} are both computed 20 times. For every run, a different realization $\bomega$ is used for the RBM-solution. The standard deviations estimated based on these 20 samples are also listed between round bracksts in Table \ref{tab:Diamond_dyn}. This table shows that the RBM-solution to \eqref{eq:RD} can be computed approximately 2 times faster as the solution to \eqref{eq:D}, and that this relative difference remains constant when the time step $h$ is varied. Table \ref{tab:Diamond_dyn} furthermore shows that the relative error in the RBM approximation reduces by a factor 2 everytime the time step $h$ is decreased by a factor 4. This is precisely the rate predicted by Theorem \ref{thm:dyn2}. 

\begin{table}[]
\centering
\begin{tabular}{|l|r|r|r|} \hline 
$\boldsymbol{h}$ & \bf{0.008} & \bf{0.002} &  \bf{0.0005} \\ \hline \hline
Sim. time (D) & 0.27 (0.03) & 1.03 (0.04) &  4.18 (0.18) \\ \hline 
Sim. time (RD) & 0.15 (0.01) & 0.57 (0.02) &  2.32 (0.06) \\ \hline 
Sim. time (RD) / Sim. time (D) [\%] & 56.81 (5.55) & 55.71 (2.52) &  55.62 (2.01) \\ \hline \hline
$|\mathbf{w}_{h,\pm}(\bomega) - \mathbf{w}_\pm|_{X} / |\mathbf{w}_\pm|_{X}$ [\%] & 7.70 (1.51) & 4.01 (0.99) &  1.99 (0.50) \\ \hline 
$|\mathbf{y}_{h}(\bomega) - \mathbf{y}|_{X} / |\mathbf{y}|_{X}$ [\%] & 4.80 (2.03) & 2.24 (1.24) &  1.11 (0.67) \\ \hline 
\end{tabular} 
    \caption{Simulation results for the randomized dynamics \eqref{eq:RD} on the diamond network in Figure \ref{fig:diamond-graph} (where $X = L^\infty(0,T; L^2(\Omega))$)}
    \label{tab:Diamond_dyn}
\end{table}

\begin{table}[]
\centering
\begin{tabular}{|l|r|r|r|} \hline 
$\boldsymbol h$ & \bf{0.008} & \bf{0.002} & \bf{0.0005} \\ \hline \hline
Sim. time (OCP) & 62.71 (0.53) & 190.47 (1.73) &  769.99 (6.32) \\ \hline 
Sim. time (ROCP) & 32.00 (4.48) & 120.50 (7.71) &  445.29 (16.60) \\ \hline 
Sim. time (ROCP) / (OCP) [\%] & 51.01 (6.93) & 63.26 (3.93) &  57.84 (2.34) \\ \hline \hline
$|J_h(\bomega,\mathbf{u}_h^*(\bomega)) - J(\mathbf{u}^*)| / J(\mathbf{u}^*)$ [\%] & 2.79 (1.43) & 1.61 (1.25) &  0.79 (0.52) \\ \hline
$|\mathbf{u}_h^*(\bomega) - \mathbf{u}|_{L^2} / |\mathbf{u}|_{L^2}$ [\%] & 22.85 (14.81) & 10.09 (5.06) &  5.35 (2.45) \\ \hline 
$|\mathbf{u}_{h}^*(\bomega) - \mathbf{u}|_{H^2} / |\mathbf{u}|_{H^2}$ [\%] & 10.88 (7.91) & 5.41 (3.33) &  2.64 (1.43) \\ \hline 
$|\mathbf{w}^*_{h,\pm}(\bomega) - \mathbf{w}^*_\pm|_{X} / |\mathbf{w}^*_\pm|_{X}$ [\%] & 23.95 (15.19) & 10.47 (5.09) &  5.61 (2.59) \\ \hline 
$|\mathbf{y}^*_{h}(\bomega) - \mathbf{y}^*|_{X} / |\mathbf{y}^*|_{X}$ [\%] & 17.85 (10.88) & 8.73 (4.47) &  4.34 (2.38) \\ \hline 
\end{tabular} 
    \caption{Simulation results for the randomized optimal control problem \eqref{eq:OCPh} on the diamond Network in Figure \ref{fig:diamond-graph} (where $X = L^\infty(0,T; L^2(\Omega))$) and $\mathbf{y}^*_h(\bomega)$ denotes the solution of \eqref{eq:D} resulting from the control $\mathbf{u}^*_h(\bomega)$)}
    \label{tab:Diamond_control}
\end{table}

In order to verify the convergence result for the optimal control in Theorem \ref{thm:contr}, we also consider an optimal control problem of the form \eqref{eq:OCP} on the same network. We choose $\alpha = 1$ and take $U_{ad} = H^2(0,T; \mathbb{R}^{|V_C|})$ and $\mathbf{y}_d(t,x) = 1$. We approximate the solution $\mathbf{u}^*(t)$ to \eqref{eq:OCP} by the solution $\mathbf{u}^*(\bomega,t)$ of the randomized optimal control problem \eqref{eq:OCPh}. The randomized dynamics are constructed in the same way as before and we again compute the solution to \eqref{eq:OCPh} 20 times in order to see the effect of the randomization. The results summarized in Table \ref{tab:Diamond_control} in which the values between round brackets again indicate the estimated standard deviations based on 20 realizations of $\bomega$. The results in Table \ref{tab:Diamond_control} show that the randomized optimal controls $\mathbf{u}^*(\bomega,t)$ are computed approximately two times faster than the optimal control $\mathbf{u}^*(t)$ to the original problem. Furthermore, the optimality gap $|J_h(\bomega,\mathbf{u}_h^*(\bomega)) - J(\mathbf{u}^*)|$, the difference $\mathbf{u}^*_h(\bomega) - \mathbf{u}^*$ in $L^2$ and $H^2$, the difference in the Riemann variables $\mathbf{w}_{h,\pm}(\bomega,t,x)$ and $\mathbf{w}_\pm(t,x)$, and the difference between the solutions $\mathbf{y}^*_h(\bomega,t,x)$ and $\mathbf{y}^*(t,x)$ of \eqref{eq:D} resulting from the controls $\mathbf{u}^*_h(\bomega,t)$ and $\mathbf{u}^*(t)$ all converge at a rate of $\sqrt{h}$, which is precisely the rate proved in Theorem \ref{thm:contr}. 

\begin{figure}
    \centering

    \subfloat[$h = 0.0080$ \label{fig:Diamond_control_00080}]{
    \includegraphics[width=0.55\textwidth]{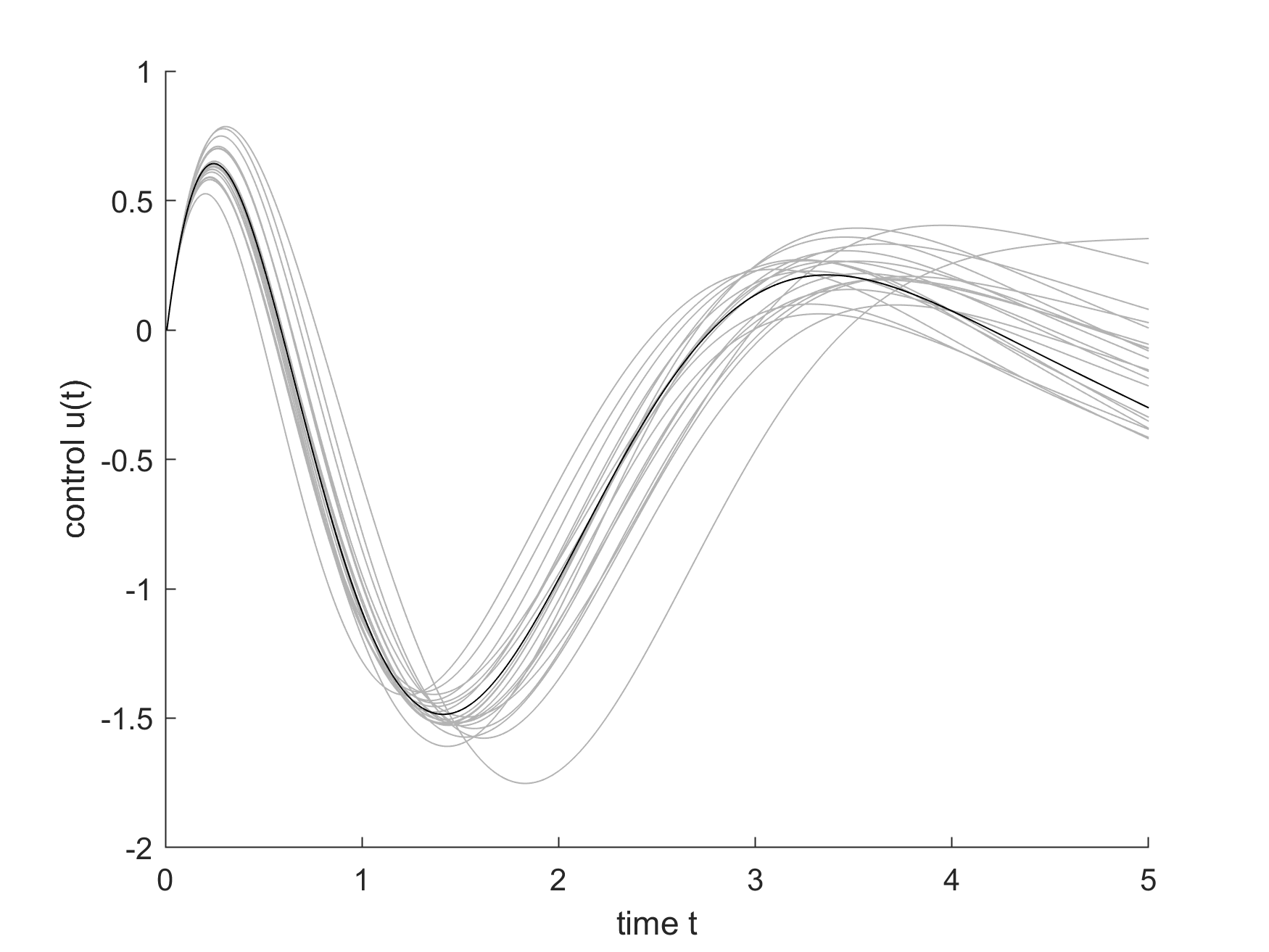}}
    \\
    \subfloat[$h = 0.0020$ \label{fig:Diamond_control_00020}]{
    \includegraphics[width=0.55\textwidth]{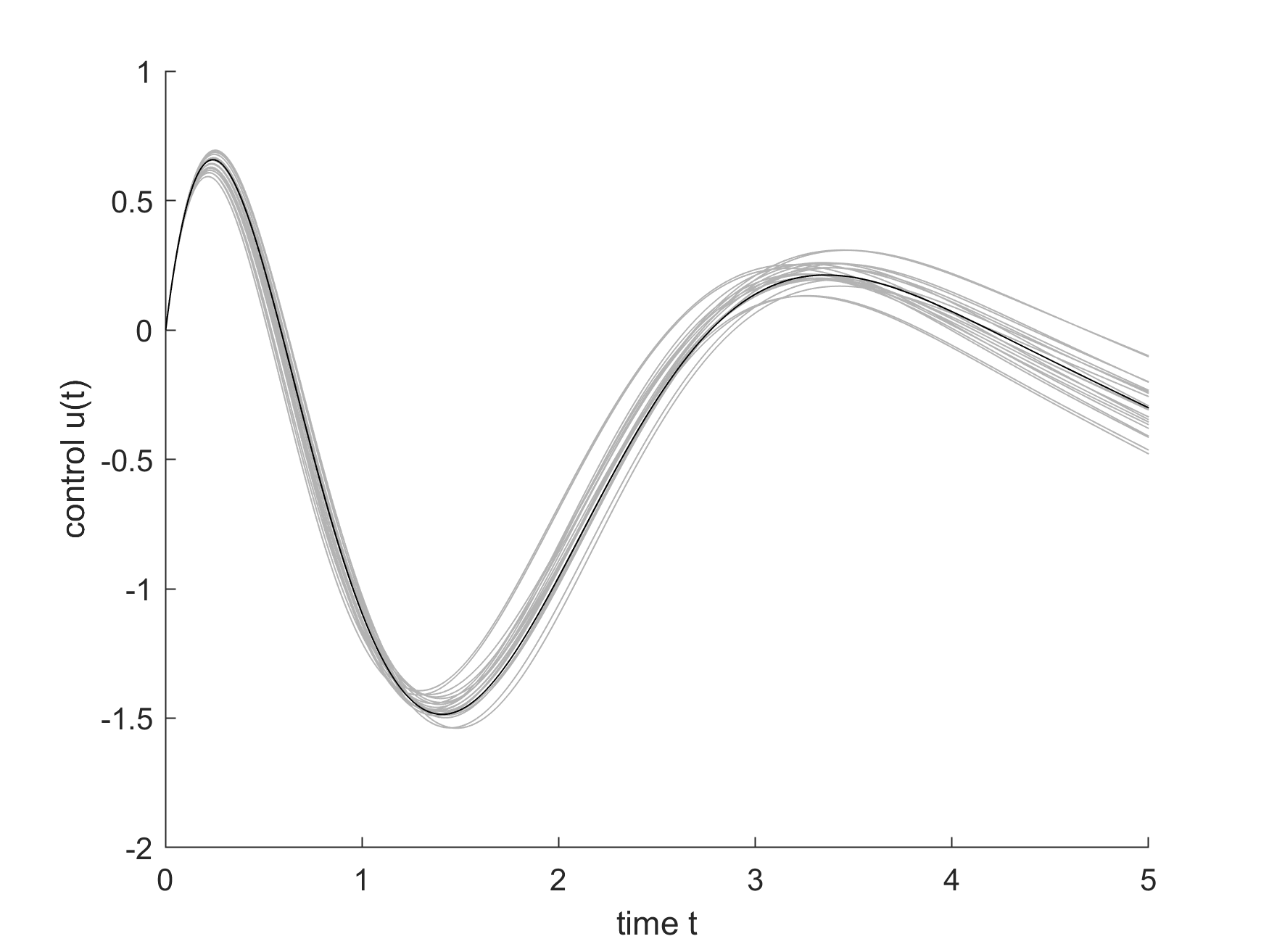}}
    \\
    \subfloat[$h = 0.0005$ \label{fig:Diamond_control_00005}]{
    \includegraphics[width=0.55\textwidth]{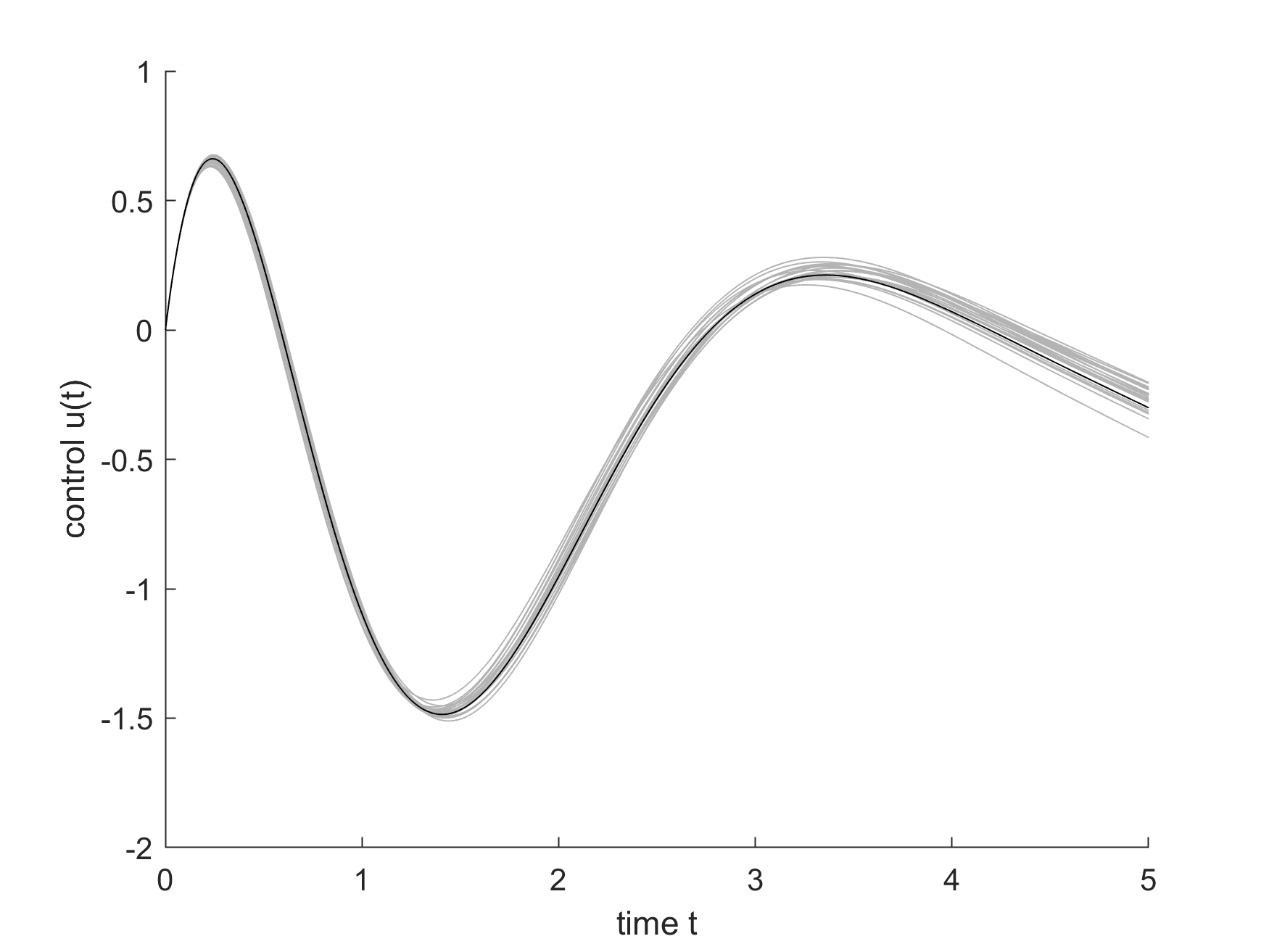}}
    \caption{The randomized optimal controls $u_h^*(\bomega,t)$ for 20 realizations of $\bomega$ (gray) compared to the optimal control $u^*(t)$ for the original problem (black) for the Diamond network from Figure \ref{fig:diamond-graph}. }
    \label{fig:Diamond_control}
\end{figure}

The optimal controls $\mathbf{u}^*_h(\bomega,t)$ and $\mathbf{u}^*(t)$ are also displayed in Figure \ref{fig:Diamond_control}. This figure also clearly shows that the optimal controls $\mathbf{u}^*_h(\bomega,t)$ for the randomized dynamics (gray) get closer and closer to optimal control $\mathbf{u}^*(t)$ for the original problem (black) when the time step $h$ is reduced. 

\subsection{Gaslib-40}
\begin{figure}
    \centering
    \includegraphics[width=\linewidth]{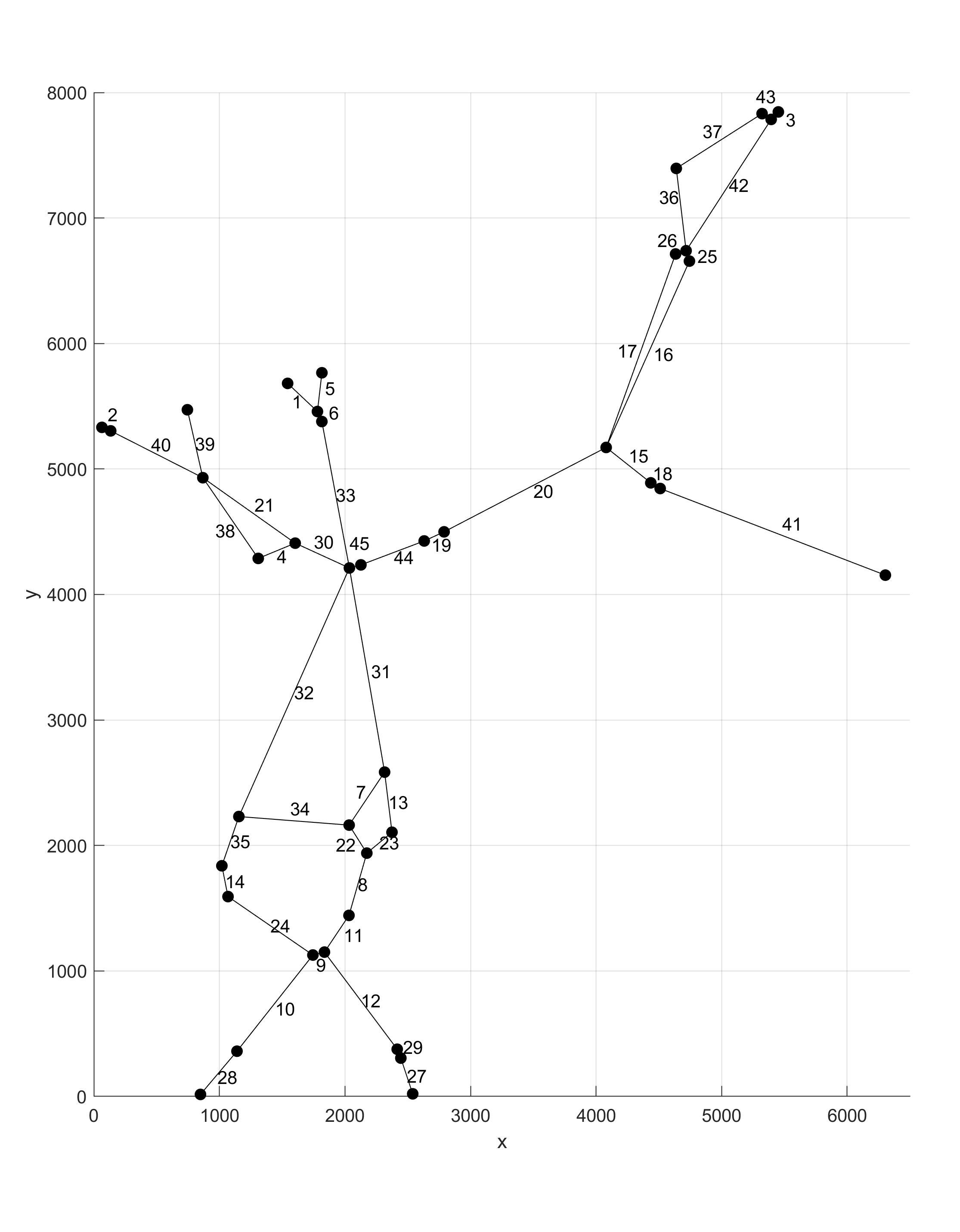}
    \caption{The network Gaslib-40 from \cite{gaslib2017, gaslibreport2017} which is simulated in the second numerical example with the numbers for each edge. 
    The network contains 10 cycles: 
    $\{ e_4,	e_{21}, e_{38} \}$, 
    $\{ e_{22},	e_{8}, e_{11}, e_{9}, e_{24}, e_{14}, e_{35}, e_{34} \}$, 
    $\{ e_{22},	e_{8}, e_{11}, e_{9}, e_{24}, e_{14}, e_{35}, e_{32}, e_{31}, e_{7} \}$, 
    $\{ e_{22},	e_{23}, e_{13}, e_{7} \}$, 
    $\{ e_{22},	e_{23}, e_{13}, e_{31}, e_{32}, e_{34} \}$, 
    $\{ e_{7},	e_{31}, e_{32}, e_{34} \}$, 
    $\{ e_{7},	e_{13}, e_{23}, e_{8}, e_{11}, e_{9}, e_{24}, e_{14}, e_{35}, e_{34} \}$, 
    $\{ e_{11},	e_{9}, e_{24}, e_{14}, e_{35}, e_{32}, e_{31}, e_{13}, e_{23}, e_{8} \}$, 
    $\{ e_{16},	e_{25}, e_{26}, e_{17} \}$, and 
    $\{ e_{36},	e_{37}, e_{43}, e_{42} \}$. 
    }
    \label{fig:gaslib40}
\end{figure}

For the second numerical example, we consider the Gaslib-40 network from \cite{gaslib2017, gaslibreport2017}, which roughly represents
a part of the German low-calorific gas transport network in the Rhine-Main-Ruhr area. The network consisting of 40 nodes and 45 edges is displayed in Figure \ref{fig:gaslib40} and contains 10 cycles. In order to improve the scaling of the problem, we reduce the $(x,y)$-coordinates in Figure \ref{fig:gaslib40} by a factor 1000. The velocity of propagation on each edge is chosen as $c_{e_i} = 5$ for all edges. 

We use the same combination as in Subsection \ref{ssec:examples_diamond} of an upwind scheme in space with backward Euler scheme in time to discretizatize \eqref{eq:D} and \eqref{eq:RD}. The spatial grid is chosen uniformly on each edge and such that the grid spacing is always below 0.02 (after rescaling the spatial coordinates by a factor 1000). For the randomized dynamics \eqref{eq:RD}, one of the 10 subsets consisting of $|E_\omega| = 26$ edges in Table \ref{tab:gaslib40active_edges} is chosen in each timt step of length $h$ with a probability $p_\omega = 1/10$. These subsets are chosen in such a way that they remove all 10 cycles present in the original network in Figure \ref{fig:gaslib40}. 

\begin{table}[]
    \centering
    \begin{tabular}{|l|} \hline
         $E_1 = \{ e_1, e_2, e_5, e_6, e_{10}, e_{11}, e_{14}, e_{15}, e_{18}, e_{19}, e_{21}, e_{22}, e_{24},$ \\
         \qquad\qquad $e_{25}, e_{27}, e_{28}, e_{30}, e_{31}, e_{33}, e_{34}, e_{37}, e_{38}, e_{40}, e_{41}, e_{43}, e_{44} \}$ \\ \hline
         $E_2 = \{ e_2, e_3, e_6, e_8, e_{10}, e_{11}, e_{13}, e_{14}, e_{16}, e_{18}, e_{20}, e_{21}, e_{23},$ \\
         \qquad\qquad $e_{25}, e_{27}, e_{28}, e_{30}, e_{31}, e_{34}, e_{35}, e_{37}, e_{38}, e_{40}, e_{41}, e_{44}, e_{45} \}$ \\ \hline
         $E_3 = \{ e_1, e_3, e_5, e_7, e_{8}, e_{10}, e_{12}, e_{14}, e_{15}, e_{17}, e_{18}, e_{20}, e_{21},$ \\
         \qquad\qquad $e_{24}, e_{25}, e_{28}, e_{29}, e_{32}, e_{33}, e_{35}, e_{37}, e_{39}, e_{40}, e_{42}, e_{43}, e_{45} \}$ \\ \hline
         $E_4 = \{ e_1, e_2, e_5, e_6, e_{8}, e_{9}, e_{11}, e_{12}, e_{14}, e_{15}, e_{17}, e_{18}, e_{20},$ \\
         \qquad\qquad $e_{21}, e_{25}, e_{27}, e_{29}, e_{30}, e_{33}, e_{34}, e_{36}, e_{37}, e_{39}, e_{40}, e_{43}, e_{44} \}$ \\ \hline
         $E_5 = \{ e_2, e_3, e_5, e_6, e_{9}, e_{10}, e_{13}, e_{14}, e_{17}, e_{18}, e_{20}, e_{21}, e_{24},$ \\
         \qquad\qquad $e_{25}, e_{27}, e_{28}, e_{30}, e_{31}, e_{33}, e_{34}, e_{37}, e_{39}, e_{41}, e_{42}, e_{44}, e_{45} \}$ \\ \hline
         $E_6 = \{ e_1, e_3, e_4, e_6, e_{8}, e_{10}, e_{11}, e_{13}, e_{14}, e_{17}, e_{18}, e_{20}, e_{21},$ \\
         \qquad\qquad $e_{24}, e_{25}, e_{27}, e_{28}, e_{30}, e_{31}, e_{33}, e_{34}, e_{37}, e_{39}, e_{41}, e_{43}, e_{45} \}$ \\ \hline
         $E_7 = \{ e_1, e_2, e_4, e_5, e_{7}, e_{9}, e_{11}, e_{12}, e_{14}, e_{15}, e_{18}, e_{19}, e_{21},$ \\
         \qquad\qquad $e_{22}, e_{25}, e_{26}, e_{28}, e_{29}, e_{32}, e_{33}, e_{35}, e_{37}, e_{40}, e_{41}, e_{43}, e_{44} \}$ \\ \hline
         $E_8 = \{ e_2, e_3, e_5, e_6, e_{8}, e_{9}, e_{12}, e_{13}, e_{15}, e_{17}, e_{19}, e_{20}, e_{24},$ \\
         \qquad\qquad $e_{25}, e_{27}, e_{28}, e_{30}, e_{31}, e_{33}, e_{35}, e_{37}, e_{38}, e_{40}, e_{41}, e_{44}, e_{45} \}$ \\ \hline
         $E_9 = \{ e_1, e_3, e_4, e_6, e_{7}, e_{9}, e_{10}, e_{12}, e_{14}, e_{16}, e_{17}, e_{19}, e_{20},$ \\
         \qquad\qquad $e_{23}, e_{24}, e_{27}, e_{28}, e_{30}, e_{31}, e_{33}, e_{35}, e_{39}, e_{40}, e_{42}, e_{43}, e_{45} \}$ \\ \hline
         $E_{10} = \{ e_1, e_2, e_4, e_5, e_{7}, e_{8}, e_{10}, e_{11}, e_{13}, e_{14}, e_{16}, e_{18}, e_{20},$ \\
         \qquad\qquad $e_{21}, e_{25}, e_{26}, e_{28}, e_{29}, e_{31}, e_{33}, e_{35}, e_{36}, e_{39}, e_{40}, e_{43}, e_{44} \}$ \\ \hline
    \end{tabular}
    \caption{The 10 subsets of active edges used in the randomized simulation of GasLib-40 Network in Figure \ref{fig:gaslib40}}
    \label{tab:gaslib40active_edges}
\end{table}

Again, the two main results in this paper are validated for this particular example. For Theorem \ref{thm:dyn2} regarding the forward dynamics, we consider the input $u(t) = \sin(4\pi t)$ applied to the controlled node $V_C = \{ v_1\}$ and start a simulation from zero initial conditions until over a time interval of length $T = 2$. The results are summarized in Table \ref{tab:Gaslib40_dyn} and show that the randomized dynamics \eqref{eq:RD} speed up the computation even more for this problem than for the diamond network. We also see that the approximation error in the Riemann invariants $\mathbf{w}_{h\pm}(\bomega,t)$ and the solution $\mathbf{y}_h(\bomega,t)$ reduces by a approximately a factor 2 every time the time step $h$ is reduced by a factor 4. This is precisely the rate of $\sqrt{h}$ we proved in Theorem \ref{thm:dyn2}. 

\begin{figure}
    \centering

    \subfloat[$h = 0.00400$ \label{fig:Network40_control_00040}]{
    \includegraphics[width=0.55\textwidth]{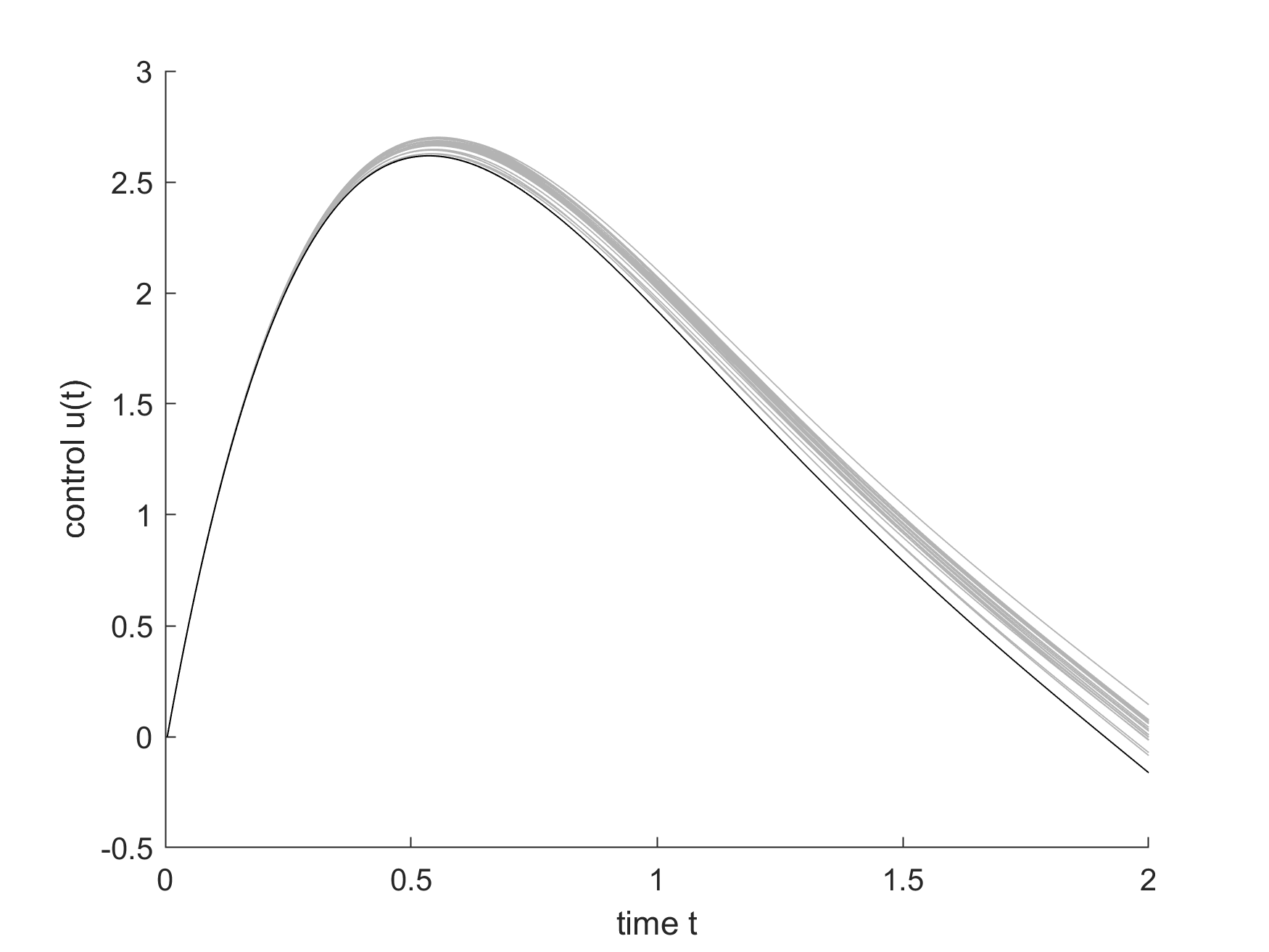}}
    \\
    \subfloat[$h = 0.00100$ \label{fig:Network40_control_00010}]{
    \includegraphics[width=0.55\textwidth]{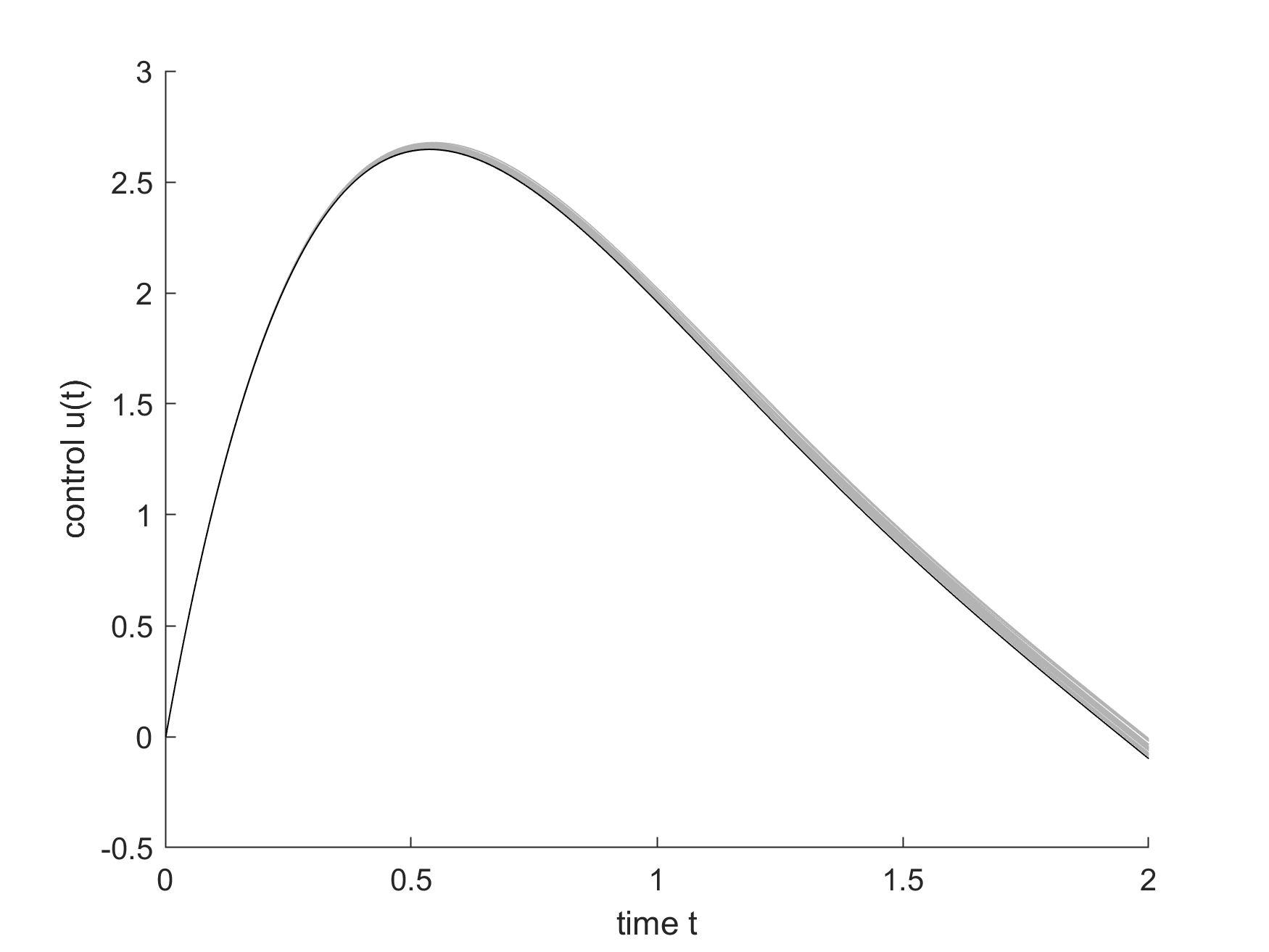}}
    \\
    \subfloat[$h = 0.00025$ \label{fig:Network40_control_000025}]{
    \includegraphics[width=0.55\textwidth]{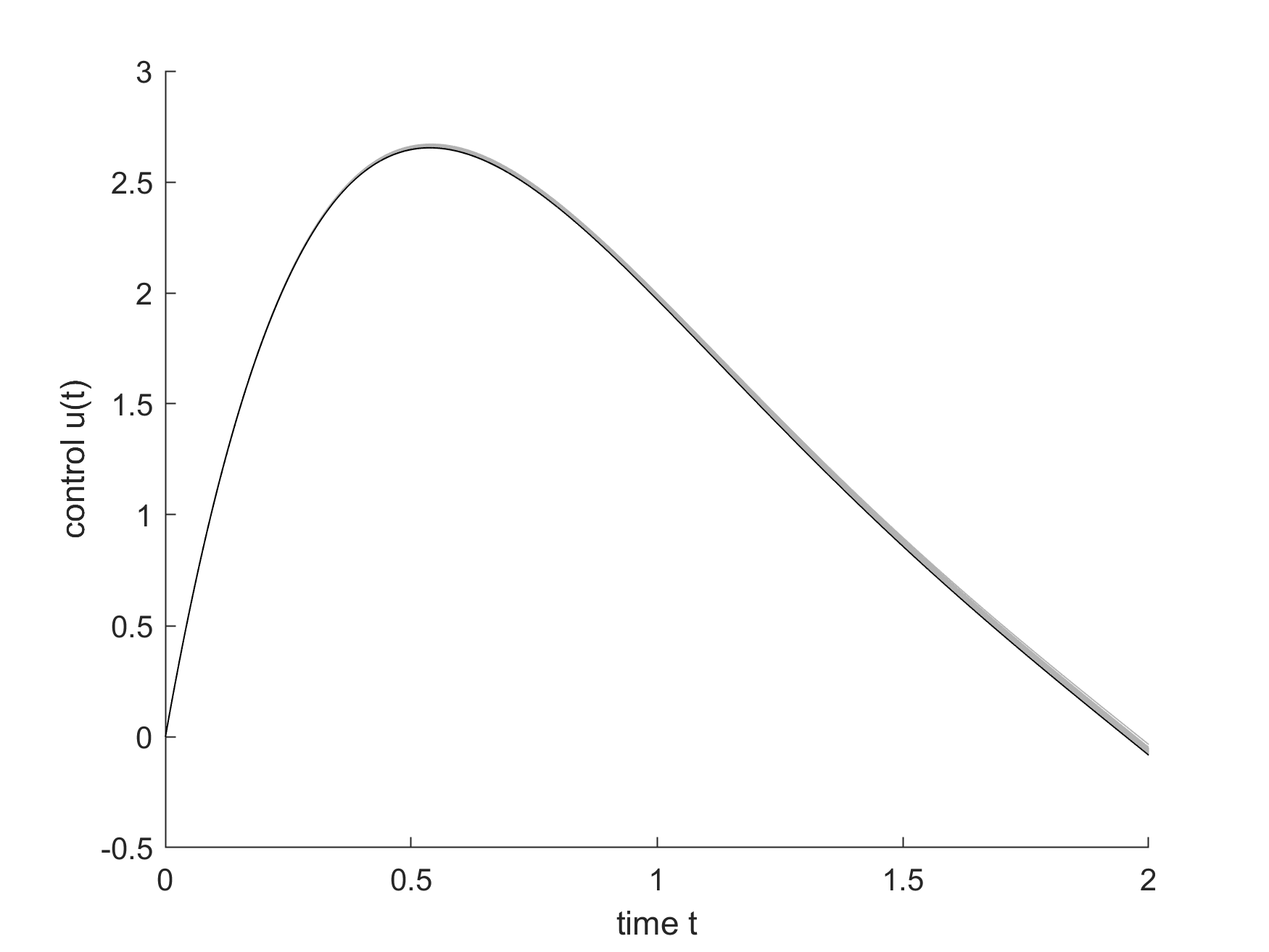}}
    \caption{The randomized optimal controls $u_h^*(\bomega,t)$ for 20 realizations of $\bomega$ (gray) compared to the optimal control $u^*(t)$ for the original problem (black) for the Gaslib-40 network displayed in Figure \ref{fig:gaslib40}. }
    \label{fig:Network40_control}
\end{figure}

We also consider also the corresponding optimal control problem \eqref{eq:OCP} on the Gaslib-40 network, again with $\alpha = 1$, $U_{ad} = H^2(0,T; \mathbb{R}^{|V_C|})$, and $\mathbf{y}_d(t,x) = 1$. We approximate the solution $\mathbf{u}^*(t)$ to \eqref{eq:OCP} again by the solution $\mathbf{u}^*_h(\bomega,t)$ of the randomized optimal control problem \eqref{eq:OCPh}. The randomization is constructed in the same way as for the forward dynamics. The results are summarized in Table \ref{tab:Gaslib40_control}. The reduction in computational time for the optimal controls is similar to the reduction observed for the forward dynamics in Table \ref{tab:Gaslib40_dyn}. We also observe that the different approximation errors related to the randomized dynamics behave like $\sqrt{h}$, which is precisely what is expected based on Theorem \ref{thm:contr}. Finally, Figure \ref{fig:Network40_control} again shows that the randomized optimal controls $\mathbf{u}^*_h(\bomega,t)$ approximate the optimal control for the original problem $\mathbf{u}^*(t)$ better and better when the time step $h$ is reduced.

\begin{table}[]
\centering
\begin{tabular}{|l|r|r|r|} \hline 
$\boldsymbol{h}$ & \bf{0.004} & \bf{0.001} & \bf{0.00025} \\ \hline \hline
Sim. time (D) & 3.32 (0.08) & 8.17 (0.18) &  32.33 (0.59) \\ \hline 
Sim. time (RD) & 1.13 (0.20) & 3.42 (0.23) &  13.37 (0.23) \\ \hline 
Sim. time (RD) / Sim. time (D) [\%] & 34.13 (6.19) & 41.86 (2.09) &  41.36 (1.06) \\ \hline \hline
$|\mathbf{w}_{h,\pm}(\bomega) - \mathbf{w}_\pm|_{X} / |\mathbf{w}_\pm|_{X}$ [\%] & 44.28 (7.13) & 26.42 (3.58) &  14.46 (2.64) \\ \hline 
$|\mathbf{y}_{h}(\bomega) - \mathbf{y}|_{X} / |\mathbf{y}|_{X}$ [\%] & 1.20 (0.31) & 0.71 (0.23) &  0.37 (0.06) \\ \hline 
\end{tabular} 
    \caption{Simulation results for the randomized dynamics \eqref{eq:RD} on the GasLib-40 network in Figure \ref{fig:gaslib40} (where $X = L^\infty(0,T; L^2(\Omega))$)}
    \label{tab:Gaslib40_dyn}
\end{table}

\begin{table}[]
\centering
\begin{tabular}{|l|r|r|r|} \hline 
$\boldsymbol h$ & \bf{0.004} & \bf{0.001} &  \bf{0.00025} \\ \hline \hline
Sim. time (OCP) & 41.49 (0.58) & 122.66 (1.02) &  376.10 (2.85) \\ \hline 
Sim. time (ROCP) & 15.30 (0.37) & 52.66 (0.80) &  183.19 (7.99) \\ \hline 
Sim. time (ROCP) / (OCP) [\%] & 36.89 (1.16) & 42.94 (0.71) &  48.71 (2.16) \\ \hline \hline
$|J_h(\bomega,\mathbf{u}_h^*(\bomega)) - J(\mathbf{u}^*)| / J(\mathbf{u}^*)$ [\%] & 0.22 (0.15) & 0.06 (0.05) &  0.04 (0.03) \\ \hline 
$|\mathbf{u}_h^*(\bomega) - \mathbf{u}|_{L^2} / |\mathbf{u}|_{L^2}$ [\%] & 7.02 (2.10) & 1.86 (0.83) &  0.73 (0.44) \\ \hline 
$|\mathbf{u}_{h}^*(\bomega) - \mathbf{u}|_{H^2} / |\mathbf{u}|_{H^2}$ [\%] & 2.21 (0.66) & 0.61 (0.25) &  0.28 (0.14) \\ \hline 
$|\mathbf{w}^*_{h,\pm}(\bomega) - \mathbf{w}^*_\pm|_{X} / |\mathbf{w}^*_\pm|_{X}$ [\%] & 5.78 (1.66) & 1.57 (0.68) &  0.59 (0.35) \\ \hline 
$|\mathbf{y}^*_{h}(\bomega) - \mathbf{y}^*|_{X} / |\mathbf{y}^*|_{X}$ [\%] & 4.47 (1.40) & 1.15 (0.54) &  0.47 (0.30) \\ \hline 
\end{tabular} 
    \caption{Simulation results for the randomized optimal control problem \eqref{eq:OCPh} on the GasLib-40 Network in Figure \ref{fig:gaslib40} (where $X = L^\infty(0,T; L^2(\Omega))$) and $\mathbf{y^*}(\bomega)$ denotes the solution of \eqref{eq:D} resulting from the control $\mathbf{u}^*_h(\bomega)$)}
    \label{tab:Gaslib40_control}
\end{table}

\section{Conclusions and Discussions \label{sec:conclusions}}
We have introduced a novel stochastic algorithm for the simulation of networks of networked wave equations, inspired by the random batch method. In this algorithm, the considered time interval is divided into subintervals of length $\leq h$. In each time interval, the dynamics are restricted to a randomly selected set of edges in the network and the other edges are frozen in time. Note that this algorithm can be interpreted as a randomized version of nonoverlapping domain decomposition. We have provided a rigorous convergence analysis for these randomized forward dynamics and the associated optimal control problems which shows that the approximation error converges to zero for $h \rightarrow 0$. In the two considered numerical examples, the computational time for both the forward dynamics and the associated optimal control problems was reduced by more than 50\% while having an approximation error of a few percent. 

Several extensions of the results in this paper are of interest. 
\begin{description}
    \item[More general hyperbolic systems] 
    Because the randomized dynamics are introduced based on the representation of the networked wave equations in terms of Riemann variables, the proofs in this paper can be extended easily to more general linear hyperbolic equations in which each Riemann variable has a constant velocity of propagation. 

    \item[Lipschitz-continuous velocity fields]
    For ease of presentation, we have restricted our analysis in this paper to the case in which the velocity of propagation fields are constant on each edge. Following the original convergence proof of the RBM in \cite{jin2020}, the convergence of the characteristics in Lemma \ref{lem:conv_xi} is also easily obtained for the case in which the velocity on each edge $c_{e_i}(x)$ is Lipschitz and decomposed into $N_{e_i}$ Lipschitz $c_{e_i,j}(x)$ as
    $$ c_{e_i}(x) = \sum_{j = 1}^{N_{e_i}} c_{e_i,j}(x).  $$
    In this way, other versions of the RBM for hyperbolic systems could be obtained that could be interpreted that as randomized overlapping domain decomposition. Note that the analysis in our manuscript is more challenging because it contains discontinuities of the velocity field at the nodes of the network which could be avoided in this Lipschitz continuous setting. 
    
    \item[Weaker regularisation]
    The optimal control problems in this paper involve a relatively strong $H^2$-regularisation. We use this regularisation to guarantee that that the optimal controls for the randomized system are Lipschitz. Whether weaker regularisation terms can be used depends on whether the Lipschitz-continuity of the optimal controls for the randomized system can also be obtained in other ways, for example from the inherent properties of the optimality system for the randomized optimal controls. This, however, is a challenging question due to the vanishing velocity of propagation (i.e., degeneration) that is inherent in the randomized dynamics. 

    \item[Nonlinear hyperbolic systems]
    The extension from general linear hyperbolic systems to semi-linear hyperbolic systems with Lipschitz nonlinearities seems straightforward. The quasi-linear case is more interesting to investigate, maybe first for specific examples such as the Euler equations or shallow water equations. 
    
    \item[Multi-dimensional spatial domains] 
    The convergence of characteristics is easily extended to more spatial dimensions. However, the treatment of the boundary conditions requires more care and would be challenging and interesting to explore further. 

    \item[Connections to Deep Learning] 
    There are strong connections between solving optimal control problems with the random batch method and the training problem for deep neural networks by stochastic gradient descent, see, e.g.,\ \cite{Geshkovski2022}. The application of the RBM to networked systems considered in this paper could be interesting to explore further in the context of models for the complex interactions between neurons in the brain, see, e.g., \cite{bassett2017}. 
\end{description}

\section*{Acknowledgments}

\backmatter

\bmhead{Acknowledgements}

We thank prof. E. Zuazua for the ideas and discussions that inspired this work. 

Daniel Veldman and Yue Wang have been supported by the Alexander von Humboldt-Professorship program awarded to E. Zuazua. Yue Wang was also supported by the Deutsche Forschungsgemeinschaft (DFG) within the Project ?Analysis and Control of Nonlinear Hyperbolic Systems with Degeneration on Networks" (no. 504042427) and within the collaborative research center TRR 154 ?Mathematical modeling, simulation and optimization using the example of gas networks" and by the Shanghai Natural Science Foundation under Grant no. 25ZR1404013. 

\begin{appendices}

\section{Networked Wave Equations}

\subsection{Derivation of BCs} \label{app:bcs}
Subtracting the two equations in \eqref{eq:riemann_inout}  yields
\begin{equation}
   2c_{e_i} D_{ij} y_x^{e_i}(t,v_j) = w_{\mathrm{out}}^{e_i}(t,v_j) - w_{\mathrm{in}}^{e_i}(t,v_j). 
\end{equation}
Therefore, the flux boundary condition can be rewritten as
\begin{equation}
   \bar{u}^{v_j}(t) = \sum_{e_i \in E(v_j)} c_{e_i}^2 D_{ij} y_x^{e_i}(t,x) = \sum_{e_i \in E(v_j)} c_{e_i} \frac{w_{\mathrm{out}}^{e_i}(t,v_j) - w_{\mathrm{in}}^{e_i}(t,v_j)}{2}.
\end{equation}
Because $w_{\mathrm{out}}^{e_i}(t,v_j) + w_{\mathrm{in}}^{e_i}(t,v_j) = w_{\mathrm{out}}(t,v_j) + w_{\mathrm{in}}(t,v_j)$  is independent of the $e_i \in E(v_j)$, 
\begin{align}
\bar{u}^{v_j}(t) &= \sum_{e_i \in E(v_j)} c_{e_i} \frac{2w_{\mathrm{out}}^{e_i}(t,v_j) - (w_{\mathrm{out}}(t,v_j) + w_{\mathrm{in}}(t,v_j))}{2} \nonumber \\
&= \sum_{e_i \in E(v_j)} c_{e_i} w_{\mathrm{out}}^{e_i}(t,v_j) -c_{\mathrm{tot},v_j} \frac{w_{\mathrm{out}}(t,v_j) + w_{\mathrm{in}}(t,v_j)}{2},
\end{align}
where $c_{\mathrm{tot},v_j} = \sum_{e_i \in E(v_j)} c_{e_i}$. Bringing $w_{\mathrm{in}}^{e_i}(t, v_j)$ and $\bar{u}^{v_j}(t)$ to the other side and dividing by $c_{\mathrm{tot},v_j} / 2$, it follows that
\begin{equation}
   w_{\mathrm{in}}^{e_i}(t,v_j) = - w_{\mathrm{out}}^{e_i}(t,v_j) + \frac{2}{c_{\mathrm{tot},v_j}}\left(- \bar{u}^{v_j}(t) + \sum_{e_i \in E(v_j)} c_{e_i} w_{\mathrm{out}}^{e_i}(t,v_j) \right), 
\end{equation}
which is the boundary condition in \eqref{eq:D}. 

\subsection{Proof of Theorem \ref{thm:Lipschitz}} \label{app:thm_Lipschitz}

The proof goes by tracing back the Riemann invariants along characteristics. 

Fix a $\delta > 0$ that is so small that characteristics cannot travel from one endpoint of an edge to the other in a time interval of length $\delta$. 
We will then prove the following induction step: if there exists a constant $\mathrm{Lip}(\mathbf{w}(k\delta))$ such that for all $e_i \in E$ and $x'', x' \in (0, \ell_{e_i})$
\begin{align}
    |w_+^{e_i}(k\delta, x'') - w_+^{e_i}(k\delta, x')| \leq \mathrm{Lip}(\mathbf{w}(k\delta)) |x'' - x'|,  \\
    |w_-^{e_i}(k\delta, x'') - w_-^{e_i}(k\delta, x')| \leq \mathrm{Lip}(\mathbf{w}(k\delta)) |x'' - x'|, 
\end{align}
then there exists constant $\mathrm{Lip}(\mathbf{w}((k+1)\delta))$ such that for all $e_i \in E$,  $x'' , x' \in  (0, \ell_{e_i})$ and $t'', t' \in [k\delta, (k+1)\delta]$
\begin{align}
    |w_+^{e_i}(k\delta, x'') - w_+^{e_i}(k\delta, x')| \leq \mathrm{Lip}(\mathbf{w}((k+1)\delta)) \left( |x'' - x'| + \bar{c} |t'' - t'| \right),  \\
    |w_-^{e_i}(k\delta, x'') - w_-^{e_i}(k\delta, x')| \leq \mathrm{Lip}(\mathbf{w}((k+1)\delta)) \left( |x'' - x'| + \bar{c} |t'' - t'| \right), 
\end{align}
where $\bar{c} = \max_{e_i \in E} c_{e_i}$. 

Fix $e_i \in E$ and choose two points $(t', x')$ and $(t'', x'')$ inside $[k\delta, (k+1)\delta] \times (0, \ell_{e_i})$. Four cases need to distinguished depending on whether each of the characteristics $s \mapsto (s, \xi^{e_i}_+(s; t'',x''))$ and $s \mapsto (s, \xi^{e_i}_+(s; t',x'))$ trace back to the initial time $\{ k\delta \} \times (0,\ell_{e_i})$ or to the boundary $[k\delta, t] \times (0, \ell_{e_i})$.  

In the first case, both characteristics trace back to $\{ k\delta \} \times (0, \ell_{e_i})$ and it follows that
\begin{equation}
w_+^{e_i}(t'',x'') - w_+^{e_i}(t', x') = w_+^{e_i}(k\delta,\xi_+^{e_i}(k\delta; t'', x'')) - w_+^{e_i}(k\delta, \xi_+^{e_i}(k\delta; t', x')). 
\end{equation}
Because $\xi^{e_i}_+(s; t,x) = x + c_{e_i}(t-s)$,
\begin{align}
|w_+^{e_i}(t'',x'') - w_+^{e_i}(t', x')|  &\leq \mathrm{Lip}(\mathbf{w}(k\delta)) \left|\xi_+^{e_i}(k\delta; t'', x'') - \xi_+^{e_i}(k\delta; t', x') \right| \nonumber \\
&\leq \mathrm{Lip}(\mathbf{w}(k\delta)) \left( |x'' - x'| + c_{e_i} |t'' - t'| \right).
\end{align}
A similar estimate holds for $w_-^{e_i}$ in the case in which the characteristics $s \mapsto (s, \xi^{e_i}_-(s; t'',x''))$ and $s \mapsto (s, \xi^{e_i}_-(s; t',x'))$ trace back to $\{ k\delta \} \times (0, \ell_{e_i})$. 

In the second case, $s \mapsto (s, \xi^{e_i}_+(s; t'',x''))$ and $s \mapsto (s, \xi^{e_i}_+(s; t',x'))$ both trace back to the boundary $[k \delta, t] \times \{ 0 \}$. For brevity, denote $t_{+,\mathrm{in}}'' :=  t^{e_i}_{+,\mathrm{in}}(t'', x'')$ and $t_{+,\mathrm{in}}' :=  t^{e_i}_{+,\mathrm{in}}(t', x')$. By using the boundary conditions in \eqref{eq:D}, it follows that
\begin{align}
&w_+^{e_i}(t'',x'') - w_+^{e_i}(t', x') 
= w_+^{e_i}(t_{+,\mathrm{in}}'',0) - w_+^{e_i}(t_{+,\mathrm{in}}',0) \nonumber \\
&= w_{\mathrm{in}}^{e_i}(t_{+,\mathrm{in}}'',v_j) - w_{\mathrm{in}}^{e_i}(t_{+,\mathrm{in}}',v_j) \nonumber \\
&= - w_{\mathrm{out}}^{e_i}(t_{+,\mathrm{in}}'',v_j) + w_{\mathrm{out}}^{e_i}(t_{+,\mathrm{in}}',v_j)  \\
& + \frac{2}{c_{\mathrm{tot}, v_j}} \left( - \bar{u}^{v_j}(t_{+,\mathrm{in}}'') + \bar{u}^{v_j}(t_{+,\mathrm{in}}') + \sum_{e_\ell \in E(v_j)} c_{e_\ell} \left( w_{\mathrm{out}}^{e_\ell}(t_{+,\mathrm{in}}'',v_j) - w_{\mathrm{out}}^{e_\ell}(t_{+,\mathrm{in}}',v_j) \right) \right). \nonumber 
\end{align}
Using that the choice of $\delta$ is such that every $w_{\mathrm{out}}^{e_\ell}$ traces back to $ \{k\delta \} \times (0, \ell_{e_i})$ which means the estimates from the first case apply also to 
$w_{\mathrm{out}}^{e_\ell}$ so that
\begin{align}
&|w_+^{e_i}(t'',x'')  - w_+^{e_i}(t', x'| \nonumber \\
&\leq \left( c_{e_i}\mathrm{Lip}(\mathbf{w}(k\delta)) + \tfrac{2}{c_{\mathrm{tot},v_j}}\left( \mathrm{Lip}(\mathbf{u}) +  \sum_{e_\ell \in E(v_j)} c_{e_\ell}^2\mathrm{Lip}(\mathbf{w}(k\delta))\right) \right) |t_{+,\mathrm{in}}'' - t_{+,\mathrm{in}}'| \nonumber \\
&\leq \left( 3\bar{c}\mathrm{Lip}(\mathbf{w}(k\delta)) + \tfrac{2}{c_{\mathrm{tot},v_j}}\left( \mathrm{Lip}(\mathbf{u})  \right) \right)|t_{+,\mathrm{in}}'' - t_{+,\mathrm{in}}'| \nonumber \\
&\leq \left( 3\bar{c}\mathrm{Lip}(\mathbf{w}(k\delta)) + \tfrac{2}{c_{\mathrm{tot},v_j}}\mathrm{Lip}(\mathbf{u}) \right) \left( |x'' - x'|/c_{e_i} + |t'' - t'| \right), \nonumber
\end{align}
where $\bar{c} = \max_{e_i \in E} c_{e_i}$ and the last step uses that $t^{e_i}_{+,\mathrm{in}}(t,x) = t - x/c_{e_i}$. 

In the third and fourth case, one of the characteristics traces back to the boundary $[k \delta, (k+1)\delta] \times \{ 0 \}$ and the other to $\{ k\delta \} \times (0, \ell_{e_i})$. Without loss of generality, consider therefore the case in which $s \mapsto (s, \xi_+(s; t'', x''))$ traces back to the boundary $[k\delta,t] \times \{ 0 \}$ and $s \mapsto (s, \xi_+(s; t', x'))$ traces back to $\{ k\delta \} \times (0, \ell_{e_i})$. In that case,
\begin{align*}
w_+^{e_i}&(t'',x'') - w_+^{e_i}(t', x') = w_+^{e_i}(t_{+,\mathrm{in}}'',0) - w_+^{e_i}(k\delta, \xi_+^{e_i}(k\delta; t', x')) \\
&=  w_{\mathrm{in}}^{e_i}(t_{+,\mathrm{in}}'',0) - w_{\mathrm{in}}^{e_i}(k\delta,0) + w_{+}^{e_i}(k\delta,0) - w_+^{e_i}(k\delta, \xi_+^{e_i}(k\delta; t', x'))
\end{align*}
Using similar estimates as in the second and first case, it follows that
\begin{align*}
|w_+^{e_i}(t'',x'') - w_+^{e_i}(t',x')| &\leq \left( 3\bar{c}\mathrm{Lip}(\mathbf{w}(k\delta)) + \tfrac{2}{c_{\mathrm{tot},v_j}}\mathrm{Lip}(\mathbf{u}) \right) |t_{+,\mathrm{in}}'' - k\delta| \\
&\qquad + \mathrm{Lip}(\mathbf{w}(k\delta)) |0 - \xi_+^{e_i}(k\delta; t', x')| \\
&\leq \left( 3\bar{c}\mathrm{Lip}(\mathbf{w}(k\delta)) + \tfrac{2}{c_{\mathrm{tot},v_j}}\mathrm{Lip}(\mathbf{u}) \right) |t_{+,\mathbf{in}}'' - t_{+,\mathbf{in}}'| 
\\
&\qquad + \mathrm{Lip}(\mathbf{w}(k\delta)) |\xi_+^{e_i}(k\delta; t'', x'') + \xi_+^{e_i}(k\delta; t', x')|,
\end{align*}
where it has been used that $t_{+,\mathrm{in}}'$ and $\xi^{e_i}_+(k\delta; t'', x'')$ are negative in this case. The remaining terms can be estimated in therms of $|x'' - x'|$ and $|t'' - t'|$ as before. 

Summarizing the results, it follows that in any of the four cases for $t \in [k\delta, (k+1)\delta]$
\begin{equation}
|w_\pm^{e_i}(t'',x'') - w_\pm^{e_i}(t', x')| \leq \left( 4\tfrac{\bar{c}}{\underline{c}}\mathrm{Lip}(\mathbf{w}(k\delta)) + \tfrac{2}{\bar{c}c_{\mathrm{tot},v_j}} \mathrm{Lip}(\mathbf{u})\right) \left( |x'' - x'|  + \bar{c} |t'' - t'| \right),
\end{equation}
where $\underline{c} = \min_{e_i \in E} c_{e_i}$. Theorem \ref{thm:Lipschitz} follows because the initial conditions in \eqref{eq:D} show that $\mathrm{Lip}(\mathbf{w}(0))= \bar{c}\mathrm{Lip}(\mathbf{y}_{0,x}) + \mathrm{Lip}(\mathbf{y}_1)$. 

\end{appendices}


\bibliography{references}

\end{document}